\documentclass[11pt]{article}
\usepackage{amsfonts}
\usepackage{amssymb}
\usepackage{graphicx}
\usepackage{mathtools}

\newtheorem{theorem}{Theorem}
 \newtheorem{corollary}[theorem]{Corollary}
 \newtheorem{lemma}[theorem]{Lemma}
 \newtheorem{proposition}[theorem]{Proposition}
\newtheorem{definition}[theorem]{Definition}
\newtheorem{remark}[theorem]{Remark}
\newtheorem{example}[theorem]{Example}

\newcommand{\Img}{\mbox{\rm im\,}}

\newcommand{\diag}{\mbox{\rm diag}}

\newcommand{\cH}{\mathcal{H}}
\newcommand{\bma}{\begin{bmatrix}}
\newcommand{\ema}{\end{bmatrix}}

\newcommand{\X}{\mathcal{X}}
\newcommand{\D}{\mathcal{D}}

\newcommand{\pperp}{\perp \!\!\!\perp}
\newcommand{\bq}{\begin{equation}}
\newcommand{\eq}{\end{equation}}
\DeclareMathOperator{\rank}{rank}

\def\BibTeX{{\rm B\kern-.05em{\sc i\kern-.025em b}\kern-.08em
    T\kern-.1667em\lower.7ex\hbox{E}\kern-.125emX}}

\newcommand{\Real}{\mathbb{R}}
\newcommand{\Comp}{\mathbb{C}}

\newcommand{\eps}{\varepsilon}

\newcommand{\proof}{\par\noindent{\bf Proof}. \ignorespaces}
\newcommand{\eproof}{\space
    {\ \vbox{\hrule\hbox{\vrule height1.3ex\hskip0.8ex\vrule}\hrule}}\par}

\newcommand {\mycomment}[1]{} 
\newcommand {\mat}  [1] {\left[\begin{array}{#1}}
\newcommand {\rix}      {\end{array}\right]}

\catcode`@=11     
\font\tenex=cmex10 
\newdimen\p@renwd
\setbox0=\hbox{\tenex B} \p@renwd=\wd0 
\def\bmat#1{\begingroup \m@th
  \setbox\z@\vbox{\def\cr{\crcr\noalign{\kern2\p@\global\let\cr\endline}}%
    \ialign{$##$\hfil\kern2\p@\kern\p@renwd&\thinspace\hfil$##$\hfil
      &&\quad\hfil$##$\hfil\crcr
      \omit\strut\hfil\crcr\noalign{\kern-\baselineskip}%
      #1\crcr\omit\strut\cr}}%
  \setbox\tw@\vbox{\unvcopy\z@\global\setbox\@ne\lastbox}%
  \setbox\tw@\hbox{\unhbox\@ne\unskip\global\setbox\@ne\lastbox}%
  \setbox\tw@\hbox{$\kern\wd\@ne\kern-\p@renwd\left[\kern-\wd\@ne
    \global\setbox\@ne\vbox{\box\@ne\kern2\p@}%
    \vcenter{\kern-\ht\@ne\unvbox\z@\kern-\baselineskip}\,\right]$}%
  \null\;\vbox{\kern\ht\@ne\box\tw@}\endgroup}
\catcode`@=12    
\newif\ifMDlatex
\catcode`@=11     
\def\MD@us#1{\csname#1style\endcsname}
\def\MD@uf#1{\csname#1font\endcsname}
\def\MD@t{text}
\def\MD@s{script}
\def\MD@ss{scriptscript}
\newdimen\MD@unit
\MD@unit\p@
\def\MD@changestyle#1{
  \relax\MD@unit0.1\fontdimen6\MD@uf{#1}0
  \everymath\expandafter{\the\everymath\MD@us{#1}}
}
\def\MD@dot{$\m@th\ldotp$}
\def\MD@palette#1{\mathchoice{#1\MD@t}{#1\MD@t}{#1\MD@s}{#1\MD@ss}}
\def\MD@ddots#1{{\MD@changestyle{#1}%
  \mkern1mu\raise7\MD@unit\vbox{\kern7\MD@unit\hbox{\MD@dot}}%
  \mkern2mu\raise4\MD@unit\hbox{\MD@dot}%
  \mkern2mu\raise \MD@unit\hbox{\MD@dot}\mkern1mu}}%
\def\MD@iddots#1{{\MD@changestyle{#1}%
  \mkern1mu\raise \MD@unit\hbox{\MD@dot}%
  \mkern2mu\raise4\MD@unit\hbox{\MD@dot}%
  \mkern2mu\raise7\MD@unit\vbox{\kern7\MD@unit\hbox{\MD@dot}}}}%
\def\MD@vdots#1{\vbox{\MD@changestyle{#1}%
    \baselineskip4\MD@unit\lineskiplimit\z@
    \kern6\MD@unit\hbox{\MD@dot}\hbox{\MD@dot}\hbox{\MD@dot}}}%
\ifMDlatex
  \DeclareRobustCommand\ddots{\mathinner{\MD@palette\MD@ddots}}%
  \DeclareRobustCommand\iddots{\mathinner{\MD@palette\MD@iddots}}%
  \DeclareRobustCommand\vdots{\mathinner{\MD@palette\MD@vdots}}%
\else
  \def\ddots{\mathinner{\MD@palette\MD@ddots}}%
  \def\iddots{\mathinner{\MD@palette\MD@iddots}}%
  \def\vdots{\mathinner{\MD@palette\MD@vdots}}%
\fi
\catcode`@=12    

\begin{document}

\title
{Differential-algebraic systems with dissipative Hamiltonian structure}
\author{
\and V. Mehrmann\footnotemark[3]~\footnotemark[1]
\and  A.J. van der Schaft \footnotemark[2]
}
\maketitle
\begin{abstract}
Different representations  of dissipative Hamiltonian and port-Hamiltonian differential-algebraic equations (DAE) systems are presented and compared. Using global geometric and algebraic points of view,  translations between the different representations are presented. Characterizations are also derived when a general DAE system can be transformed into one of these structured representations. Approaches for computing the structural information and the described transformations are derived that can be directly implemented as numerical methods. The results are demonstrated with a large number of examples.
\end{abstract}
\noindent

\noindent


\renewcommand{\thefootnote}{\fnsymbol{footnote}}
\footnotetext[1]{
Institut f\"ur Mathematik, Sekr.~MA 4-5, TU Berlin, Stra{\ss}e des 17.~Juni 136,
10623 Berlin, Germany.
\texttt{mehrmann@math.tu-berlin.de}. Research supported by Deutsche Forschungsgemeinschaft DFG through priority Programm SPP 1984 within  project 790 39-2: Distributed dynamic control of network security.
}
\footnotetext[2]{Bernoulli institute for Mathematics, Computer Science and AI, Jan C. Willems Center for Systems and Control, University of Groningen, Nijenborgh 9, Groningen, the Netherlands}
\renewcommand{\thefootnote}{\arabic{footnote}}

\section{Introduction}\label{sec:intro}
Since the original introduction of  the energy based modeling concept of port-Hamiltonian (pH) systems in \cite{MasS92,SchM95}, see also \cite{Bre08,GolSBM03,JacZ12,OrtSMM01,Sch00,Sch06,SchJ14}, various novel definitions and formulations have been made to incorporate on the one hand systems defined on manifolds, see \cite{EbeMS07,SchM18} and on the other hand systems with algebraic constraints defined in the form of Differential-Algebraic Equations (DAEs), see \cite{BeaMXZ18,MehM19}.
Constraints in pH systems typically arise in electrical or transport networks, where Kirchhoff's laws constrain the models at network nodes, as balance equations in chemical engineering, or as holonomic or non-holonomic constraints in mechanical multibody systems. Furthermore, they arise in the interconnection of pH systems when the interface conditions are explicitly formulated and enforced via Lagrange multipliers,
see \cite{BeaMXZ18,MehM19,Sch13,SchJ14,SchM18,SchM20} for a variety of examples, or the recent survey \cite{MehU23}.

In this paper we recall  different model representations of pHDAEs, respectively dissipative Hamiltonian DAEs (dHDAEs), and provide a systematic geometric and algebraic theory that also provides a 'translation' between the different representations. 
We also discuss approaches for explicitly computing the structural information that can be directly implemented as numerical methods. We mainly restrict ourselves to finite-dimensional linear time-invariant systems without inputs and outputs, but indicate in several places where extensions to general systems are possible.  We also discuss only real systems although most of the results can also be formulated for complex problems in a similar way.

For general linear time-invariant homogeneous differential-algebraic equations
{
 \begin{equation}
\frac{d}{dt} (Ex)= E\dot x= Ax, \label{gendae}
\end{equation}
where $E,A\in\mathbb R^{n,n}$ (without additional geometric and algebraic  structures), as well as for their extensions to linear time-varying and nonlinear systems, the theory is well understood \cite{KunM06}.
On the other hand, in contrast to general DAEs \eqref{gendae}, the DAEs that arise in energy based modeling have extra structure and symmetries, and it is natural to identify and exploit this extra structure for purposes of analysis, simulation and control.

In this paper we will restrict our attention to \emph{regular}  systems, i.e., systems with
$\det (\lambda E-A)$ not identically zero, although many of the results are expected to extend to over- and underdetermined systems; cf. the Conclusions.
The structural properties for regular systems are characterized via the \emph{(real) Weierstra{\ss}}  canonical form of the matrix pair $(E,A)$, see e.g. \cite{Gan59a}. To determine this canonical form one computes nonsingular matrices $U,W$ such that
\[
UEW=\mat{cc} I & 0 \\ 0 & \mathcal N\rix,\ UAW= \mat{cc} \mathcal J & 0 \\ 0 & I \rix,
\]
where $\mathcal J$ is in real Jordan canonical form and $\mathcal N$ is nilpotent in Jordan canonical form which is associated to the eigenvalue $\infty$

An important quantity that will be used throughout the paper is the size of the largest Jordan block in $\mathcal N$ (the index of nilpotency), which is called the \emph{(differentiation)  index of the pair $(E,A)$ as well as the associated DAE}, where, by convention,  $\nu=0$ if $E$ is invertible.
}
{

The goal of this paper is to study in detail the relationships and differences between different classes of (extended) Hamiltonian DAEs; without making any a priori invertibility assumptions. Furthermore, we exploit at the same time a \emph{geometric} point of view on Hamiltonian DAEs using Lagrange and Dirac structures (as well maximally monotone subspaces; cf. Section \ref{sec:dissipation}), and an \emph{algebraic} point of view using (computationally feasible) condensed forms. From a DAE analysis point of view, similar to the result known already for Hamiltonian DAEs of the form \eqref{gendhdae} in \cite{MehMW18}, we also show that the index of an (extended) dissipative Hamiltonian DAEs can be at most two. Moreover, we will show that index two algebraic constraints may only arise from singularity of $P$, i.e., from the Lagrange structure, and thus are strictly linked to singular Hamiltonians.

Another important question that we will answer is the characterization when a general DAE is equivalent to a structured DAE the form (\eqref{gendhdae} or (\eqref{gendila}).

The paper is organized as follows.

In Section~\ref{sec:dhdaes} we discuss dissipative Hamiltonian DAEs and present several examples. The concept of extended Hamiltonian DAEs is discussed in Section~\ref{sec:eHDAEs}.
In Section~\ref{sec:not:prelim} we will present the geometric theory of dissipative Hamiltonian DAEs, we introduce Dirac and Lagrange structures as well maximally monotone subspaces to incorporate dissipation. Section~\ref{sec:geometric} presents different coordinate representations and their relation. To derive  algebraic characterizations, in Section~\ref{sec:equi} we present different types of equivalence transformations and corresponding condensed forms for the different classes of (extended dissipative) Hamiltonian DAEs. In Section~\ref{sec:daetoph} it is analyzed when general linear DAEs can be represented as  (extended dissipative) Hamiltonian systems. In all sections we present examples that illustrate the  properties of the different representations. In several appendices we present proofs that can be implemented as numerically stable algorithms.

}

\section{Dissipative Hamiltonian DAEs}\label{sec:dhdaes}
{A well established  representation of DAE systems with symmetry structure is that of \emph{dissipative Hamiltonian DAEs (dHDAEs)} which has been studied in detail in \cite{BeaMXZ18,MehMW18,MehMW21}. These systems have the form}
 \begin{equation}
{ \frac{d}{dt}(Ez)}= E\dot z= (J-R)Q z, \label{gendhdae}
\end{equation}
where  $J=-J^\top , R=R^\top  \in \mathbb R^{\ell,\ell}$, and $E,Q\in\mathbb R^{\ell,n}$ with $E^\top Q=Q^\top E$. In this representation the associated Hamilton function (\emph{Hamiltonian}) is given by a quadratic  form
\begin{equation}\label{defhzEQ}
\mathcal H^z(z)=\frac 12 z^\top E^\top Q z,
\end{equation}
which (as an energy) is typically nonnegative for all $z$, but more general Hamiltonians also arise in practice. Note that in the ODE case, i.e. if  $E=I$, then this reduces to the standard quadratic Hamiltonian $\mathcal H^z(z)=\frac 12 z^\top Q z$.

\begin{remark}\label{rem:Edotx}{\rm The equivalent formulations $\frac{d}{dt}(Ez)= E\dot z$ are both used in the literature and obviously lead to the same results in the linear constant coefficient case, but not anymore in the linear time-varying or nonlinear case. Actually, for the structured DAEs considered in this paper $E$ is the product $E=KP$ of two matrices $K$ and $P$, and for generalizations it is appropriate to consider
$K\frac{d}{dt}(Pz)$.
}
\end{remark}

Furthermore, the dissipation term $R$ is typically positive semidefinite, denoted as $R\geq 0$, but we will first discuss the case that $R=0$ which we call \emph{Hamiltonian DAE (HDAE)}.

\begin{remark}\label{rem:port1}{\rm The definition of (dHDAEs) in \eqref{gendhdae} can be easily extended to systems with ports in which case it takes the form
\begin{eqnarray}\nonumber
{ \frac{d}{dt}(Ez)}=E\dot z=&=& (J-R)Q z+ (F-G) u,\\
y&=& (F+G)^\top Q z + D u, \label{gendphDAE}
\end{eqnarray}
with
\[
\mat{cc} Q^\top J Q  & Q^\top F \\ -F^\top Q  &\frac 12 (D-D^\top)\rix=-\mat{cc} Q^\top J Q  & Q^\top F \\ -F^\top Q & \frac 12 (D-D^\top)\rix^\top
\]
and
\[
\mat{cc} Q^\top R Q  & Q^\top G \\ G^\top Q & \frac 12 (D+D^\top)\rix=\mat{cc} Q^\top R Q  & Q^\top G \\ G^\top Q & \frac 12 (D+D^\top)\rix^\top \geq 0,
\]
see \cite{BeaMXZ18}. We will call these systems \emph{dissipative pHDAE (dpHDAE) systems}.
}
\end{remark}

{

\begin{remark}\label{rem:defham}{\rm
Note that by the symmetry of $E^\top Q$ the Hamiltonian $\mathcal H^z(z)=\frac 12 z^\top E^\top Q z$ can be also written as a function of $Qz$ or $Ez$.
}
\end{remark}

\begin{remark}\label{rem:nonsymmR}{\rm
In standard port-Hamiltonian modeling the matrix $R$ in \eqref{gendhdae} is always assumed to be {\it symmetric} (as well as positive semi-definite). Alternatively, one could start from a general matrix $R$ only satisfying $R+R^\top \geq 0$. Then by splitting $R$ into its symmetric and skew-symmetric part, one could add the skew-symmetric part to the skew-symmetric structure matrix $J$ and continue as in \eqref{gendhdae} with the symmetric part $\frac 12 (R+R^\top)$ of $R$.}

\end{remark}
}

Structured DAEs of the form~\eqref{gendhdae} arise naturally in all physical domains.
\begin{example}\label{ex:circuit}{\rm 
The subclass of RLC networks as considered in e.g. \cite{BeaMXZ18,Dai89,Fre11} has the form
\begin{eqnarray}\label{eq:RLC_network_1}
\underbrace{\mat{ccc} D_C C D_C^\top & 0 & 0\\ 0 & L & 0 \\ 0 & 0 & 0 \rix}_{:=E}
\mat{c}\dot{V} \\ \dot{I_L}\\ \dot{I_S} \rix=
\underbrace{\mat{ccc} -D_R G D_R^\top & -D_L & -D_S \\ D_L^\top & 0 & 0\\ D_S^\top & 0 & 0 \rix}_{:=J-R}
\mat{c}V \\ I_L\\ I_S \rix,
\end{eqnarray}
where the real positive definite diagonal matrices $L$, $C$, $G$  describe inductances,  capacitances, and conductances, respectively. The matrices $D_C, D_L, D_R, D_S$ are the parts of the incidence matrix corresponding, respectively, to the capacitors, inductors, resistors (conductors), and current sources of the circuit graph, where $D_S$ is of full column rank. Furthermore $V$ are the node potentials and $I_L, I_S$ denote the currents through the inductors and sources, respectively.  This system has the form \eqref{gendhdae}, with $E=E^\top \geq 0$, $Q$ equal to the identity matrix, and where $J$ and $-R$ are defined to be the skew-symmetric and symmetric part,
respectively, of the matrix on the right hand side of~\eqref{eq:RLC_network_1}.
The Hamiltonian is given by $\mathcal H(V,I_L)= \frac 12 V^\top  D_C C D_C^\top V+  \frac 12 I_L^\top  L I_L$ and it does not involve the variables $I_S$, which are in the kernel of $E$.}
\end{example}

\begin{example} \label{ex:stokes} {\rm Space discretization of the Stokes equation in fluid dynamics, see, e.g., \cite{EmmM13}, leads to a dissipative Hamiltonian system
\[
\mat{cc}M&0\\0 &0\rix \mat{c} \dot v_h \\ \dot p_h \rix = \mat{cc} A & B\\ -B^\top& 0\rix \mat{c} v_h \\ p_h \rix +\mat{c} f_h \\ 0 \rix,
\]
where $A=A^\top $ is a positive semidefinite discretization of the negative Laplace operator, $B$ is a discretized gradient, and $M=M^\top$ is a positive definite mass matrix. The homogeneous system has the form \eqref{gendhdae}, with
\[
R=\mat{cc}A&0\\ 0&0\rix,\ J=\mat{cc}0&B\\ -B^\top&0\rix,\ E= \mat{cc}M&0\\0 &0\rix,\ Q=I.
\]
The Hamiltonian is given by $\mathcal H= \frac 12 v^\top _h M v_h$; it does not involve the variables $p_h$.}
\end{example}

\begin{example} \label{ex:gas}{\rm Space discretization of the Euler equation describing the acoustic wave propagation in a gas pipeline network \cite{EggK18,EggKLMM18} leads to a DAE
\[ \mat{ccc} M_1& 0 &0\\ 0 & M_2 & 0\\ 0 & 0 & 0 \rix \mat{c} \dot p_h\\ \dot q_h \\ \dot \lambda \rix=
\left ( \mat{ccc} 0& -G &0\\ G^\top  & 0 & K^\top\\ 0 & -K & 0\rix -\mat{ccc} 0& 0 &0\\ 0 & D & 0\\ 0 & 0 & 0\rix \right ) \mat{c} p_h \\ q_h \\ \lambda \rix +f,
\]
where $p_h$ is the discretized pressure, $q_h$ is a discretized flux, and $\lambda$ is a Lagrange multiplier
that penalizes the violation of the conservation of mass and momentum at the pipeline nodes.

The homogeneous system has the form \eqref{gendhdae} with $Q=I$, and Hamiltonian
$ \mathcal H(p_h,q_h) =\frac 12 ( p_h^\top M_1 p_h+ q_h^\top M_2 q_h)$. The Hamiltonian does not involve the Lagrange multiplier $\lambda$.
}
\end{example}

\begin{example}\label{ex:mech}{\rm Consider a linear mechanical system (with $q$ denoting the vector of position coordinates) $M\ddot q +D\dot q + Wq= f$, together with kinematic constraints $G\dot q=0$, see e.g. \cite{EicF98}. {The constraints give rise to {\it constraint forces} $G^\top \lambda$, with $\lambda$  a vector of Lagrange multipliers. This yields the dynamics} $M\ddot q +D\dot q + Wq  = f+G^\top  \lambda$, and the resulting system can be written in first order form as
\[
 \mat{ccc}M & 0 & 0\\ 0 & I &0 \\ 0 & 0 & 0 \rix\mat{c}\ddot q \\ \dot q \\ \dot  \lambda\rix
 +\mat{ccc}D & I & -G^\top \\ -I& 0 & 0 \\  G & 0 & 0 \rix \mat{ccc} I & 0 & 0 \\ 0 & W & 0 \\ 0  & 0 & I\rix \mat{c}\dot q\\ q \\ \lambda \rix=\mat{c}f\\0 \\ 0\rix.
\]
Here, $E=E^\top$  and $Q=Q^\top$ are commuting matrices, $E$, $W$, and $R=R^\top$ are positive semidefinite, so the homogeneous system is of the form \eqref{gendhdae}. The Hamiltonian is given by $\frac 12 ( \dot q^\top M \dot q + q^\top W q)$ (kinetic plus potential energy) and does not involve the Lagrange multiplier $\lambda$.
}
\end{example}

{ As indicated in the presented examples, singularity of $E$, and thus the presence of algebraic constraints, implies that the Hamiltonian (the total stored energy) $\mathcal H^z(z)=\frac 12 z^\top E^\top Q z$ does not involve all of the variables contained in the vector $z$.} In fact in some of the examples the variables that do not show up in the Hamiltonian are \emph{Lagrange multipliers}. Conversely, singularity of $E$ may arise as a limiting situation of otherwise regular Hamiltonians, as shown by the following simple example.

\begin{example} \label{exmsd}{\rm
{ Consider the model of standard mass-spring-damper system
with model equation
\[
\mat{c} \dot{q} \\ \dot{p} \rix = \mat{cc} 0 & 1 \\ -1 & -d \rix \mat{cc} k & 0 \\ 0 & \frac{1}{m} \rix \mat{c} q \\ p \rix
\]
and Hamiltonian $\mathcal H(q,p) = \frac{1}{2} kq^2 + \frac{p^2}{2m}$.
\begin{center}
\includegraphics[scale=1]{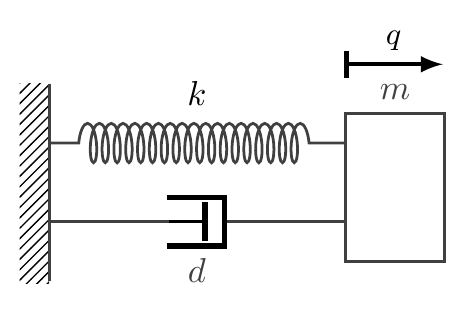}
\end{center}
To compute the limit $m \to 0$,} we rewrite the system in coordinates $q$ and $v:=\frac{p}{m}$ as
\[
\mat{c} \dot{q} \\ m\dot{v} \rix = \mat{cc} 0 & 1 \\ -1 & -d \rix \mat{c} kq \\ v \rix
\]
with Hamiltonian $\mathcal H(q,v)= \frac{1}{2} kq^2 + \frac{1}{2}mv^2$.
For $m \to 0$ this converges to the DAE system
\[
\mat{c} \dot{q} \\ 0 \rix = \mat{cc} 0 & 1 \\ -1 & -d \rix \mat{c} kq \\ v \rix,
\]
which is of (differentiation) index one if $d \neq 0$, and of index two  if $d=0$. The limiting Hamiltonian $\mathcal H(q,v)= \frac{1}{2} kq^2$ is not a function of $v$ anymore.

Alternatively, we can compute the limit for $k \to \infty$. For this we rewrite the system in coordinates $F:=kq$ and $p$ as
\[
\mat{c} \frac{1}{k}\dot{F} \\ \dot{p} \rix = \mat{cc} 0 & 1 \\ -1 & -d \rix \mat{c} F \\ \frac{p}{m} \rix,
\]
with Hamiltonian $\mathcal H(F,p)= \frac{1}{2k} F^2 + \frac{1}{2m}p^2$.
For $k \to \infty$ this converges to the DAE system
\[
\mat{c} 0 \\ \dot{p} \rix = \mat{cc} 0 & 1 \\ -1 & -d \rix \mat{c} F \\ \frac{p}{m} \rix,
\]
which is index two for any $d$, with Hamiltonian $\mathcal H(F,p)= \frac{1}{2m}p^2$, which is not involving $F$.

We may also take the limits $m \to 0$ and $k \to \infty$ simultaneously. By rewriting the system  in the variables $F,v$, with Hamiltonian $\mathcal H(F,p)= \frac{1}{2k} F^2 + \frac{1}{2}mv^2$. This leads to the purely algebraic system
\[
\mat{c} 0 \\ 0 \rix = \mat{cc} 0 & 1 \\ -1 & -d \rix \mat{c} F \\ v \rix,
\]
with zero Hamiltonian, and having a single solution $v=0$, $F=0$, irrespective of the damper.
}
\end{example}


\section{Extended Hamiltonian DAEs}\label{sec:eHDAEs}
Classically, {Hamiltonian systems} and also their extension to systems with inputs and outputs, called {port-Hamiltonian systems}, see e.g. \cite{DuiMSB09,GolSBM03,MasS92,SchM95,SchM13,Sch00,SchJ14}, are defined by a \emph{Dirac structure}, representing the power-conserving interconnection structure, an \emph{ energy-dissipation relation}, and a \emph{Hamilton function (Hamiltonian)}, capturing the total energy storage in the system. The Dirac structure formalizes the \emph{generalized junction structure} known from port-based modeling theory. Importantly, the Dirac structure may entail linear constraints on the so-called effort variables, called \emph{effort constraints}. Through the gradient vector of the Hamilton function, these effort constraints induce algebraic constraints on the state variables, see especially \cite{MehM19,Sch13} for a discussion of general nonlinear port-Hamiltonian DAEs. In the special case of a linear autonomous system without energy dissipation and without inputs and outputs, { and after choosing a basis such that $\mathcal X=\mathbb R^{n}$,}  the equations of motion take the DAE form
\begin{equation}
\label{KL}
K\dot{x}=LQx,
\end{equation}
where the pair of matrices $K,L\in \mathbb R^{n,n}$, satisfying $KL^\top  + LK^\top =0$ and $\rank \mat{cc} K& L \rix= n$, describes the Dirac structure, and $\mathcal H^x(x) = \frac{1}{2}x^\top  Q x$ defines the Hamilton function for some $Q=Q^\top $. Denoting by $\ker$ and $\Img$ the kernel and image of a linear map or its matrix representation, geometrically the Dirac structure is defined by the subspace
$\mathcal D \subset {\mathcal X} \times \mathcal X^*$ given as
$\ker \mat{cc} K & L \rix$, where $\mathcal X$ is the linear state space, and $\mathcal X^*$ its dual space.

\begin{remark}{\rm
More precisely $\mathcal D \subset \hat{\mathcal X} \times \mathcal X^*$, where $\hat{\mathcal X}$ is the \emph{tangent space} to $\mathcal X$ at $x \in \mathcal X$. However, since $\mathcal X$ is linear, tangent spaces at any $x\in \mathcal X$ can be identified with each other and with $\mathcal X$.
}
\end{remark}

{
\begin{remark}\label{rem:coord-free}{\rm
In this paper we frequently switch between coordinate free representations with state space $\mathcal X$ and coordinate representations that are obtained for the case
$\mathcal X=\mathbb R^{n}$ after choosing a basis. Furthermore, in this case it will be tacitly assumed that the dual basis is chosen for the dual space $\mathcal X^*$. Although this is not the most general setup (from an abstract linear algebraic or functional analytic point of view), it makes the presentation of results much more convenient.
}
\end{remark}
}
Algebraic constraints occur if the matrix $K$ is singular, and are represented via
\[
Qx \in \Img K^\top.
\]
Singularity of $K$ often results from the network structure as the following example demonstrates.
\begin{example} \label{ex:LCnetwork} {\rm Consider a general linear $LC$-electrical circuit. Let the circuit graph be determined by an incidence matrix $D$, defining Kirchhoff's current laws $I \in \ker D$ and voltage laws $V \in \Img D^\top $. Split the currents $I$ into currents $I_C$ through the capacitors and currents $I_L$ through the inductors, and the voltages $V$ into voltages $V_C$ across the capacitors and voltages $V_L$ across the inductors. Furthermore, let $I_C=-\dot{q}_C$ and $V_L= - \dot{\varphi}$, with $q$ the vector of charges at the capacitors and $\varphi$ the flux linkages of the inductors. Split the incidence matrix $D$ accordingly as $D= \mat{cc} D_C & D_L \rix$. Then Kirchhoff's current laws take the form
\[ 
D_C\dot{q}=D_L I_L.
\] 
Furthermore, let $F$ be a maximal annihilator of $D^\top $, i.e., $\ker F= \Img D^\top $. Then Kirchhoff's voltage laws are given as $F \mat{c} V_C \\ V_L \rix= 0$, and after splitting $F= \mat{cc} F_C & F_L \rix$ accordingly, we have
\[
F_L\dot{\varphi}=F_C V_C.
\]
Writing the linear constitutive equations for the capacitors as $q=C V_C$ for some positive definite diagonal capacitance matrix $C$ and those for the inductors as $\varphi=L I_L$ for some positive definite diagonal inductance matrix $L$, we finally obtain the system of equations
\[
\mat{cc} D_C & 0 \\ 0 & F_L \rix \mat{c} \dot{q} \\ \dot{\varphi} \rix = \mat{cc} 0 & D_L \\ F_C & 0 \rix \mat{cc} C^{-1} & 0 \\ 0 & L^{-1} \rix \mat{c} q \\  \varphi \rix,
\]
which is in the form \eqref{KL} with Hamilton function $\mathcal H(q,\varphi)= \frac{1}{2}q^\top  C^{-1} q +  \frac{1}{2}\varphi^\top  L^{-1} \varphi$. It is easily checked that singularity of $K =\mat{cc} D_C & 0 \\ 0 & F_L \rix$ corresponds to parallel interconnection of capacitors or series interconnection of inductors. (Note that since $\ker F= \Img D^\top $ the linear space $\Img F^\top = \ker D$ is spanned by the cycles of the circuit graph.)
See e.g.  \cite{GunBJR21} for pHDAE modeling of electrical circuits.
}
\end{example}
Motivated by \cite{BeaMXZ18} the port-Hamiltonian point of view on DAE systems was extended in
\cite{SchM18} by replacing the gradient vector $Qx$ of the Hamiltonian function $\cH(x)=\frac{1}{2} x^\top Qx$ by a general Lagrangian subspace, called  \emph{Lagrange structure} in the present paper, since we want to emphasize the similarity with Dirac structures. As we will see, this extension allows to bridge the gap with the (dissipative) Hamiltonian formulation \eqref{gendphDAE} of the previous Section \ref{sec:dhdaes}.

In the linear homogeneous case without dissipation, one starts with a Dirac structure $\mathcal D \subset {\mathcal X} \times \mathcal X^*$ and a Lagrange structure $\mathcal L \subset \mathcal X \times \mathcal X^*$.
The \emph{composition} of the Dirac structure $\D$ and the Lagrange structure $\mathcal L$, over the \emph{shared variables} $e \in \mathcal X^*$, is then defined as
\begin{equation}\label{composition}
\mathcal D \circ \mathcal L = \{(f,x) \in {\mathcal X} \times \mathcal X \mid \mbox{there exists } e \in \mathcal X^* \mbox{ such that } (f,e) \in \mathcal D,\ (e,x) \in \mathcal L \}.
\end{equation}
Substituting $f=-\dot{x}$ this leads to the coordinate-free definition of the dynamics
\begin{equation}
\label{pH}
(-\dot{x},x) \in \mathcal D \circ \mathcal L.
\end{equation}
In order to obtain a coordinate representation of the dynamics \eqref{pH}, the simplest option (later on in Section \ref{sec:effortflow} we will discuss another one) is to start from the \emph{image representation} of the Lagrange structure $\mathcal L$, defined by a pair of matrices $P,S\in \mathbb R^{n,n}$ satisfying $P^\top S=S^\top P$ and $\rank \mat{cc} P^\top  & S^\top  \rix= n$. Taking coordinates $x$ for $\mathcal X=\mathbb R^{n}$ and dual coordinates $e$ for its dual space $\mathcal X^*=\mathbb R^{n}$, the image representation of $\mathcal L$ is given by
\begin{equation}
\label{imageL}
\mat{c} x \\ e \rix = \mat{c} Pz \\  Sz \rix,
\end{equation}
for some parameterizing vector $z$ in a space $\mathcal Z$ (of the same dimension as $\mathcal X$). Analogously, we consider a \emph{kernel representation} of the Dirac structure $\D$ given by matrices $K,L$ satisfying $KL^\top=-LK^\top$ and $\rank \mat{cc} K  & L  \rix= n$, such that $\D = \{(f,e) \in \mathcal X \times \mathcal X^* \mid Kf + Le=0 \}$. (More details will be given in Section \ref{sec:not:prelim}.)

Substituting $f=-\dot{x} = P\dot{z}$ and $e=Sz$ this leads to the coordinate representation
\begin{equation}
KP\dot z= LS z. \label{gendila}
\end{equation}
We will call this class \emph{extended Hamiltonian differential-algebraic systems (extended HDAEs)}. If dissipation is incorporated, see Section~\ref{sec:dissipation}, then it is called \emph{extended dHDAEs}.

Importantly, the presence of algebraic constraints in \eqref{gendila} may arise \emph{both} by singularity of $K$ (as was already the case for \eqref{KL}) as well as by singularity of $P$ (as was the case for \eqref{gendhdae}). This motivated the introduction of the notions of \emph{Dirac algebraic constraints} (corresponding to singularity of $K$) and of \emph{Lagrange algebraic constraints} (corresponding to singularity of $P$) in \cite{SchM18}.

Similar to dHDAE systems \eqref{gendhdae}, the Hamiltonian of the extended HDAE system \eqref{gendila} is specified by the Lagrange structure, and is given by
\begin{equation}\label{defHzPS}
\mathcal H^z(z)= \frac 12 z^\top S^\top Pz.
\end{equation}
Indeed, one immediately has the \emph{energy conservation property}
\[
\frac{d}{dt} \mathcal H^z(z) = z^\top S^\top P\dot{z} = - e^\top f =0,
\]
since $(f,e) \in \D$ {and thus $e^\top f=0$}.

While singularity of $K$ in physical systems modeling typically arises from \emph{interconnection} due to the network structure, singularity of $P$ often arises as a limiting situation. An elaborate example will be provided later as Example \ref{ex2mass}.
\begin{remark}\label{rem:firsttrafo}{\rm
Note that if $K$ is \emph{invertible}, then by multiplying \eqref{gendila} with $K^{-1}$ from the left, we obtain the lossless version of the dHDAE system \eqref{gendhdae} in Section \ref{sec:dhdaes} with $E=P$, $Q=S$, $J=-J^\top = K^{-1} L$ (and $R=0$). Thus the replacement of the Hamiltonian $\cH(x)=\frac 12 x^\top  Q x$ by a Lagrange subspace \eqref{imageL} constitutes a first step towards an \emph{overarching} formulation of (dissipative) Hamiltonian DAE systems. 
}
\end{remark}

\section{Geometric theory of (dissipative) Hamiltonian DAEs}\label{sec:not:prelim}
In this section we take a systematic geometric view on Hamiltonian DAE systems, extending the existing geometric treatment of extended HDAE systems, as already discussed  in Section \ref{sec:eHDAEs}. We also incorporate the discussed classes of dHDAE systems from  Section \ref{sec:dhdaes}.

Consider an n-dimensional  linear state space $\mathcal X$ with elements denoted by $x$. Let $\hat{\mathcal X}$ denote the \emph{tangent space} to $\mathcal X$ at $x \in \mathcal X$, with elements denoted by $f$ and called \emph{flow} vectors.
{As mentioned before}, since $\X$ is linear, tangent spaces at different $x \in \mathcal X$ can be identified with each other and with $\mathcal X$; implying that $f \in \X$ as well.
Furthermore, let $\mathcal X^*$ be the \emph{dual  space} of $\mathcal X$, with elements denoted by $e$ and called \emph{effort} vectors.

\subsection{Dirac and Lagrange structures}\label{sec:dila}
The product space $\mathcal X \times \mathcal X^*$ is endowed with the two canonical bilinear forms
\begin{equation}\label{bilin}
\begin{array}{l}
\langle (f_1,e_1), (f_2,e_2) \rangle_+ := f_1^\top e_2 + f_2^\top e_1, \\[2mm]
\langle (f_1,e_1), (f_2,e_2) \rangle_- := f_1^\top e_2 - f_2^\top e_1,
\end{array}
\end{equation}
represented by the two matrices
\[
\Pi_+ := \mat{cc} 0 & I_n \\ I_n & 0 \rix, \quad \Pi_- :=\mat{cc} 0 & I_n \\ -I_n & 0 \rix,
\]
where we recognize $\Pi_-$ as the standard symplectic form on $\mathcal X$.

\begin{definition}\label{def:dirlag}
A subspace $\mathcal D \subset \mathcal X \times \mathcal X^*$ is called a 
\emph{Dirac structure} if the bilinear form $\langle \cdot, \cdot \rangle_+$ is zero on $\mathcal D$ and moreover $\mathcal D$ is maximal with respect to this property. A subspace $\mathcal L \subset \mathcal X \times \mathcal X^*$ is a \emph{Lagrange structure} if the bilinear form $\langle \cdot, \cdot \rangle_-$ is zero on $\mathcal L$ and moreover $\mathcal L$ is maximal with respect to this property. A Lagrange structure $\mathcal L \subset \mathcal X \times \mathcal X^*$ is called \emph{nonnegative} if the quadratic form defined by $\Pi_+$ is nonnegative on $\mathcal L$.
\end{definition}
\begin{remark} \label{rem:drilag} {\rm
In this paper we have chosen the terminology 'Lagrange structure', instead of the more common terminology 'Lagrangian subspace', in order to emphasize the similarity to Dirac structures. Also note that the definition of Dirac structures can be extended to manifolds instead of linear state spaces $\mathcal X$; in this context Dirac structures on linear spaces are often referred to as \emph{constant} Dirac structures.}
\end{remark}

We have the following characterizations of Lagrange and Dirac structures, see e.g. \cite{DuiMSB09,SchJ14}.
\begin{proposition}\label{prop:ladi}
Consider an $n$-dimensional  linear state space $\mathcal X$ and its dual space $\mathcal X^*$.
\begin{itemize}
\item [i)] A subspace $\mathcal D \subset \mathcal X \times \mathcal X^*$ is a Dirac structure if and only if $\mathcal D=\mathcal D^{\pperp_+}$, where ${}^{\pperp_+}$ denotes the orthogonal complement with respect to the bilinear form $\langle \cdot, \cdot \rangle_+$. Furthermore, $\mathcal D \subset \mathcal X \times \mathcal X^*$ is a Dirac structure if and only $e^\top f=0$ for all $(f,e) \in \mathcal D$ and $\dim \mathcal D=n$.

\item [ii)] A subspace $\mathcal L \subset \mathcal X \times \mathcal X^*$ is a Lagrange structure if and only if $\mathcal L=\mathcal L^{\pperp_-}$, where ${}^{\pperp_-}$ denotes the orthogonal complement with respect to the bilinear form $\langle \cdot, \cdot \rangle_-$. Any Lagrange structure satisfies $\dim \mathcal L =n$.
\end{itemize}
\end{proposition}

Dirac and Lagrange structures admit structured \emph{coordinate representations}, see e.g. \cite{SchJ14,SchM18}. For  this paper the following representations are most relevant.
Using matrices $K,L \in \mathbb R^{n,n}$, any \emph{Dirac structure} $\mathcal D \subset {\mathcal X} \times \mathcal X^*$ admits the kernel/image representation
\begin{equation}\label{diracspace}
\mathcal D = \ker \mat{cc} K & L \rix = \Img \mat{cc} L^\top \\ K^\top \rix \subset {\mathcal X} \times \mathcal X^*,
\end{equation}
with $K,L$ satisfying $\rank \mat{cc} K & L\rix=n$ and the generalized \emph{skew-symmetry condition}
\begin{equation}\label{dirac}
 KL^\top + L K^\top =0.
\end{equation}
Conversely any such pair $K,L$ defines a Dirac structure.

Analogously, any Lagrange structure $\mathcal L \subset \mathcal X \times \mathcal X^*$ can be represented as
\begin{equation}\label{lagrangespace}
\mathcal L = \ker \mat{cc} -P^\top & S^\top \rix = \Img \mat{cc} S \\ P \rix \subset \mathcal X^* \times \mathcal X,
\end{equation}
for certain matrices $S,P  \in \mathbb R^{n,n}$ satisfying $\rank \mat{c} P \\ S \rix =n$ as well as the generalized \emph{symmetry condition}
\begin{equation} \label{lagrange}
S^\top P= P^\top S.
\end{equation}
A Lagrange structure is, furthermore, \emph{nonnegative} if and only if $S^\top P \geq 0$.

As already described in Section~\ref{sec:eHDAEs}, by using the image representation $x=Pz$, $e=Sz$ of the Lagrange structure $\mathcal L$, and the kernel representation $Kf + Le=0$ of the Dirac structure $\mathcal D$ one is led to the representation \eqref{gendila} of the extended HDAE system defined by $\mathcal D$ and $\mathcal L$.

The following is a physical example where both $K$ and $P$ turn out to be singular. The singularity of $K$ is due to the presence of kinematic constraints, while the singularity of $P$ is caused by a limiting argument in the energy expression.

\begin{example}\label{ex2mass}{\rm

Consider two masses $m_1$ and $m_2$ connected by a spring with spring constant $k$, where the right mass $m_2$ is subject to the kinematic constraint $v_2=0$ (velocity is zero). With positions $q_1,q_2$ and momenta $p_1,p_2$, the Hamiltonian is given by
\[
\cH(q_1,q_2,p_1,p_2)= \frac{1}{2}k(q_1-q_2)^2 + \frac{p_1^2}{2m_1} + \frac{p_2^2}{2m_2}.
\]
Denoting by $e_{q1}, e_{q2}$ the spring forces at both ends of the spring, and  by $e_{p1}, e_{p2}$ the velocities of the two masses, we obtain the relation
\[
\begin{bmatrix} q_1 \\q_2 \\ p_1 \\p_2 \\\hline e_{q1} \\ e_{q2} \\  e_{p1} \\ e_{p2} \end{bmatrix} =
\begin{bmatrix}
1 & 0 & 0 & 0 \\
0 & 1 & 0 & 0 \\
0 & 0 & 1 & 0 \\
0 & 0 & 0 & 1 \\
\hline
k & -k & 0 & 0 \\
-k & k & 0 & 0 \\
0 & 0 & \frac{1}{m_1} & 0 \\
0 & 0 & 0 & \frac{1}{m_2} \\
\end{bmatrix}
\begin{bmatrix} z_1 \\ z_2 \\ z_3 \\ z_4 \end{bmatrix}
\]
To consider the limit $k \to \infty$, meaning that the spring is replaced by rigid connection, we first express the system in different coordinates.
\[
\begin{bmatrix}
z_1\\ z_2 \\z_3 \\z_4
\end{bmatrix}
=
\begin{bmatrix}
\frac{m_2}{m_1 + m_2} & 1 & 0 & 0 \\
- \frac{m_1}{m_1 + m_2} & 1 & 0 & 0 \\
0 & 0 & 1 &\frac{m_1}{m_1 + m_2}  \\
0 & 0 & -1 &\frac{m_2}{m_1 + m_2}
\end{bmatrix}
\begin{bmatrix}
\tilde{z}_1\\ \tilde{z}_2 \\ \tilde{z}_3 \\ \tilde{z}_4
\end{bmatrix}.
\]
This yields the transformed representation
\[
\begin{bmatrix} q_1 \\q_2 \\ p_1 \\p_2 \\ \hline e_{q1} \\ e_{q2} \\  e_{p1} \\ e_{p2} \end{bmatrix} =
\begin{bmatrix}
\frac{m_2}{m_1 + m_2} & 1 & 0 & 0 \\
-\frac{m_1}{m_1 + m_2} & 1 & 0 & 0 \\
0 & 0 & 1 & \frac{m_1}{m_1 + m_2} \\
0 & 0 & -1 & \frac{m_2}{m_1 + m_2} \\
\hline
k & 0 & 0 & 0 \\
-k & 0 & 0 & 0 \\
0 & 0 & \frac{1}{m_1} & \frac{1}{m_1+m_2}  \\
0 & 0 & -\frac{1}{m_2} & \frac{1}{m_1 + m_2}
\end{bmatrix}
\begin{bmatrix} \tilde{z}_1 \\ \tilde{z}_2 \\ \tilde{z}_3 \\ \tilde{z}_4 \end{bmatrix}.
\]
Taking the limit $k \to \infty$ yields the Lagrange structure $\mathcal L$ in image representation
\[
\begin{bmatrix} q_1 \\q_2 \\ p_1 \\p_2 \\ \hline e_{q1} \\ e_{q2} \\  e_{p1} \\ e_{p2} \end{bmatrix} =
\begin{bmatrix}
0 & 1 & 0 & 0 \\
0 & 1 & 0 & 0 \\
0 & 0 & 1 & \frac{m_1}{m_1 + m_2} \\
0 & 0 & -1 & \frac{m_2}{m_1 + m_2} \\
\hline
1 & 0 & 0 & 0 \\
-1 & 0 & 0 & 0 \\
0 & 0 & \frac{1}{m_1} & \frac{1}{m_1+m_2}  \\
0 & 0 & -\frac{1}{m_2} & \frac{1}{m_1 + m_2}
\end{bmatrix}
\begin{bmatrix} \tilde{z}_1 \\ \tilde{z}_2 \\ \tilde{z}_3 \\ \tilde{z}_4 \end{bmatrix}
\]
and the limiting Hamiltonian is just the kinetic energy
\[
\cH(\tilde{z}_1,\tilde{z}_2,\tilde{z}_3,\tilde{z}_4) = \frac{1}{2}\left(\frac{1}{m_1} + \frac{1}{m_2}\right)\tilde{z}^2_3 + \frac{\tilde{z}^2_4}{2(m_1 + m_2)}=\frac{\tilde{z}^2_4}{2(m_1 + m_2)}.
\]
The system has a \emph{Lagrange algebraic constraint} due to the linear dependency in the rows of~$P$.
The Dirac structure $\mathcal D$ is given as
\begin{eqnarray*}
\mathcal D &=& \left \{(f_{q1}, f_{q2}, f_{p1}, f_{p2}, e_{q1}, e_{q2}, e_{p1}, e_{p2} =0) \mid  \mbox{there exists }\lambda \right . \\
&& \mbox{ s.t. }
\begin{bmatrix} f_{q1}\\ f_{q2}\\ f_{p1}\\ f_{p2} \end{bmatrix} =
\begin{bmatrix} 0 & 0& -1 & 0 \\  0 & 0& 0 & -1 \\  1 & 0& 0 & 0 \\  0 & 1 & 0 & 0 \end{bmatrix}
\begin{bmatrix} e_{q1}\\ e_{q2}\\ e_{p1}\\ e_{p2} \end{bmatrix} +
\begin{bmatrix} 0 \\ 0 \\ 0 \\ 1 \end{bmatrix} \left .\lambda \right \}.
\end{eqnarray*}
After elimination of the Lagrange multiplier $\lambda$ this yields
\[
K= \begin{bmatrix} 1 & 0 & 0 & 0 \\ 0 & 1 & 0 & 0 \\ 0 & 0 & 1 & 0 \\ 0 & 0 & 0 & 0 \end{bmatrix}
, \quad
L= \begin{bmatrix} 0 & 0 & 1 & 0 \\ 0 & 0 & 0 & 1 \\ -1 & 0 & 0 & 0 \\ 0 & 0 & 0 & 0 \end{bmatrix},
\]
and hence
\[
KP = \begin{bmatrix} 0 & 1 & 0 & 0 \\ 0 & 1 & 0 & 0 \\ 0 & 0 & 1 & \frac{m_1}{m_1 + m_2} \\ 0 & 0 & 0 & 0 \end{bmatrix},
\quad
LS = \begin{bmatrix} 0 & 0 & \frac{1}{m_1} & \frac{1}{m_1 + m_2} \\ 0 & 0 & -\frac{1}{m_2} & \frac{1}{m_1 + m_2} \\ -1 & 0 & 0 & 0 \\ 0 & 0 & - \frac{1}{m_2} & \frac{1}{m_1 + m_2} \end{bmatrix}.
\]
Finally subtracting the second equation from the first equation, we  obtain the HDAE system
\[
\begin{bmatrix}
0 & 0 & 0 & 0 \\
0 & 1 & 0 & 0 \\
0 & 0 & 1 & \frac{m_1}{m_1 + m_2} \\
0 & 0 & 0 & 0
\end{bmatrix}
\begin{bmatrix} \dot{\tilde{z}}_1 \\ \dot{\tilde{z}}_2 \\ \dot{\tilde{z}}_3 \\ \dot{\tilde{z}}_4 \end{bmatrix} =
\begin{bmatrix}
0 & 0 & \frac{1}{m_1} + \frac{1}{m_2} & 0 \\
0 & 0 &  - \frac{1}{m_2} & \frac{1}{m_1 + m_2}\\
-1 & 0 & 0 & 0 \\
0 & 0 & - \frac{1}{m_2} & \frac{1}{m_1 + m_2}
\end{bmatrix}
\begin{bmatrix} \tilde{z}_1 \\ \tilde{z}_2 \\ \tilde{z}_3 \\ \tilde{z}_4 \end{bmatrix}.
\]
Here the first equation is the Lagrange algebraic constraint $\tilde{z}_3=0$ (and eventually $\tilde{z}_1=0$) obtained by letting $k \to \infty$ (corresponding to singularity of $P$), and the last equation is the Dirac algebraic constraint $- \frac{\tilde{z}_3}{m_2} + \frac{\tilde{z}_4}{m_1 + m_2}=0$, i.e., $\tilde{z}_4=0$ resulting from the kinematic constraint, leading to singularity of $K$ and resulting in the trivial dynamics $\dot{\tilde{z}}_2  =  0 $.
}
\end{example}

\begin{remark}\label{rem:FErepresent}{\rm
Instead of using a parametrization of the Lagrange structure $\mathcal L$ one can also use a parameterization of the Dirac structure $\mathcal D$,
\begin{equation}\label{diracz}
\mat{c} f \\ e \rix= \mat{c} L^\top \\ K^\top \rix v
\end{equation}
with $v \in \mathcal V$, where $\mathcal V$ is an $n$-dimensional parameter space. This yields an extended dHDAE system (but now in the parameter vector $v$) given by
\begin{equation} \label{gendiladual}
P^\top K^\top \dot{v}= S^\top L^\top v,
\end{equation}
which is the \emph{adjoint system} of \eqref{gendila}. See \cite{KunM11} for a detailed discussion of adjoint systems of DAEs.
}
\end{remark}

\subsection{Incorporation of dissipation}\label{sec:dissipation}
As noted in Section \ref{sec:not:prelim}, extended HDAE systems \eqref{gendila}, geometrically defined by a Dirac and Lagrange structure, already include HDAE systems \eqref{gendhdae} \emph{without} dissipation. Conversely, any HDAE system with $K$ invertible can be rewritten into the form \eqref{gendhdae} with $R=0$.

In order to complete the geometric viewpoint towards the inclusion of dissipation (and thus to \eqref{gendhdae}), we recall the geometric definition of a port-Hamiltonian system \cite{SchM95,Sch17,SchJ14}. By replacing the Hamiltonian function by a Lagrange structure as in \cite{SchM18}, and specializing to the case without external variables (inputs and outputs), such systems will be called \emph{extended dHDAE systems}.
\begin{definition}\label{def:edhadae}
Consider a state space $\mathcal X$ with linear coordinates $x$ and a linear space of resistive flows $\mathcal F_R$. Furthermore, consider a Dirac structure $\mathcal D$ on $\mathcal X \times \mathcal F_R$, a Lagrange structure $\mathcal L \subset \mathcal X \times \mathcal X^*$, and a nonnegative Lagrange structure $\mathcal R \subset \mathcal F_R \times \mathcal F^*_R$. Then an \emph{extended dissipative Hamiltonian DAE (extended dHDAE) system} is defined as the  tuple $(\mathcal X,\mathcal F_R,\mathcal D,\mathcal L,\mathcal R)$ with
\begin{equation}\label{edhdae}
\begin{array}{rcl}
&&\left \{ (\dot{x},x)  \mid  \mbox{there exist } e \in \mathcal X^*, f_R \in \mathcal F_R, e_R \in \mathcal F^*_R \right . \\[2mm]
&& \quad \left . \mbox{ such that } (-\dot{x},-f_R,e,e_R) \in \mathcal D, (f_R,e_R) \in \mathcal R, \, (x,e) \in \mathcal L \right \}.
\end{array}
\end{equation}
\end{definition}
If $\mathcal L$ is represented as in \eqref{lagrangespace}, i.e.
\begin{equation}
\label{SP1}
\mathcal L= \left \{\mat{c} x \\ e \rix \mid \mat{c} x \\ e \rix = \mat{c} P \\ S \rix z, \; z \in \mathbb R^n \right \},
\end{equation}
then it immediately follows from the properties of the Dirac structure $\mathcal D$ and the nonnegative Lagrange structure $\mathcal R$ that the dynamics of the extended dHDAE satisfies
\[
\frac{d}{dt} (\frac{1}{2} z^\top S^\top P z) = e_R^\top f_R \leq 0.
\]
More generally we will now introduce the notion a \emph{maximally monotone subspace}, which is overarching the notions of a Dirac structure $\mathcal D$ and a nonnegative Lagrange structure $\mathcal R$.
%
\begin{definition}\label{def:maxmono}
Consider a linear space $\mathcal X$. A subspace $\mathcal M \subset \mathcal X \times \mathcal X^*$ is  called \emph{monotone} subspace if
\begin{equation}\label{mono}
e^\top f \geq 0
\end{equation}
for all $(f,e) \in \mathcal M$, and it is \emph{maximally monotone} if additionally $\mathcal{M}$ is maximal with respect to this property (i.e., there does not exist a monotone subspace $\mathcal M' \subset \mathcal X \times \mathcal X^*$ with $\mathcal M \subsetneq \mathcal M'$).
\end{definition}
\begin{remark} \label{rem:monrel} {\rm The definition of a monotone subspace is a special case of the notion of a monotone \emph{relation} $\widetilde{M}$, which is defined as a \emph{subset} of  $\mathcal X \times \mathcal X^*$ satisfying
\begin{equation} \label{mono1}
(e_1 - e_2)^\top (f_1 - f_2) \geq 0
\end{equation}
for all $(f_1,e_1), (f_2,e_2) \in \widetilde{M}$. Clearly if $\widetilde{M}$ is a subspace then \eqref{mono1} reduces to \eqref{mono}. (Maximally) monotone subspaces with a sign change were employed before in \cite{GerHR21}, using the terminology of '(maximally) linear dissipative relations'. Nonlinear port-Hamiltonian systems with respect to a general (maximally) monotone relation were coined as \emph{incremental port-Hamiltonian systems} in \cite{CamS13}; see \cite{CamS22} for further developments.
}
\end{remark}

Obviously, a subspace $\mathcal M$ is monotone if and only if the quadratic form defined by $\Pi_+$ is nonnegative on $\mathcal M$, since
$\langle (f,e), (f,e) \rangle_+ = 2 e^\top f$. This yields
\begin{proposition}\label{prop:repmono}
Consider a state space $\mathcal X$ with  $\dim \mathcal X=n$. Then any monotone subspace of $\mathcal X \times \mathcal X^*$ has dimension less than or equal to $n$, and any maximally monotone subspace of $\mathcal X \times \mathcal X^*$ has dimension $n$. Any maximally monotone subspace $\mathcal M \subset \mathcal X \times \mathcal X^*$ can be represented as
\begin{equation} \label{MN}
\mathcal M = \Img \mat{c} N^\top \\ M^\top \rix
\end{equation}
for $M,N\in \mathbb R^{n,n}$ satisfying $\rank \mat{cc} N & M \rix=n$ and
\begin{equation} \label{MN1}
MN^\top + NM^\top \geq 0.
\end{equation}
Conversely, any subspace defined by $M,N$ satisfying \eqref{MN1} is a maximally monotone subspace.
\end{proposition}
\proof
The proof follows, since $\Pi_+$ has $n$ positive and $n$ negative eigenvalues.
\eproof

Obviously any Dirac structure $\mathcal D$ given by a pair of matrices $K,L$ is maximally monotone (take $M=K$ and $N=L$). In a similar way any nonnegative Lagrange structure $\mathcal R$ given by a pair of matrices $P,S$ with $S^\top P \geq 0$ is maximally monotone by taking $N^\top =P$, $M^\top =S$.

Importantly, also the \emph{composition} of two maximally monotone subspaces is again maximally monotone.
In order to prove this we first state the following lemma.
\begin{lemma}\label{lem:pullback}
Let $A: \mathcal F \to \mathcal G$ be a linear map between two linear spaces $\mathcal F,\mathcal G$. Let $\mathcal M_{\mathcal G} \subset \mathcal G \times \mathcal G^*$ be a maximally monotone subspace. Then the \emph{pull-back of $\mathcal M_{\mathcal G}$ via $A$}, defined as
\[
b_A (\mathcal M_{\mathcal G}) := \{ (f,A^*g^*) \in \mathcal F \times \mathcal F^* \mid (Af,g^*) \in \mathcal M_{\mathcal G} \},
\]
is maximally monotone. Furthermore, let $\mathcal M_{\mathcal F} \subset \mathcal F \times \mathcal F^*$ be a maximally monotone subspace. Then the \emph{push-forward of $\mathcal M_{\mathcal F}$ via $A$}, defined as
\[
f_A (\mathcal M_{\mathcal F}) := \{ (Af,g^*) \in \mathcal G \times \mathcal G^* \mid (f,A^*g^*) \in \mathcal M_{\mathcal F} \},
\]
is maximally monotone.
\end{lemma}
\proof
It is immediately checked that $b_A (\mathcal M_{\mathcal G})$ is monotone. Furthermore
\[
\begin{array}{l}
\dim b_A (\mathcal M_{\mathcal G}) = \dim \mathcal M_{\mathcal G} + \dim \ker A - \dim \ker A^* = \\[2mm]
\dim \mathcal G + \dim \ker A - \dim \ker A^* = \dim \mathcal F,
\end{array}
\]
and thus $b_A (\mathcal M_{\mathcal G})$ is maximally monotone. The proof to show that  $f_A (\mathcal M_{\mathcal F})$ is maximally monotone is analogous.
\eproof

Using Lemma~\ref{lem:pullback} we can show that maximally monotone subspaces satisfy the following composition property. This same property was recently derived for maximally monotone {\it relations} in \cite{CamS22}, assuming additional regularity conditions.
\begin{proposition}
\label{prop:comp} Consider an extended dHDAE system as in \eqref{edhdae} and let
$\mathcal M_a$ and $\mathcal M_b$ be maximally monotone subspaces
\[
\mathcal M_a \subset  \mathcal F_a \times \mathcal F \times \mathcal E_a \times \mathcal E, \quad \mathcal M_b \subset  \mathcal F_b \times \mathcal F \times \mathcal E_b \times \mathcal E
\]
with $\mathcal E=\mathcal F^*, \mathcal E_a=\mathcal F_a^*, \mathcal E_b=\mathcal F_b^*$.
Define the composition
\[
\begin{array}{rcl}
\mathcal M_a \circ \mathcal M_b & := &  \{(f_a,f_b,e_a,e_b) \mid \mbox{ there exist } f \in \mathcal F, e \in \mathcal E \mbox{ such that } \\[2mm]
&& (f_a,f,e_a,e) \in \mathcal M_a, (f_b,-f,e_b,e) \in \mathcal M_b \}.
\end{array}
\]
Then $\mathcal M_a \circ \mathcal M_b \subset \mathcal F_a \times \mathcal F_b \times \mathcal E_a \times \mathcal E_b$ is again maximally monotone.
\end{proposition}
\proof
Let $\mathcal V_a:=\mathcal F$ and $\mathcal V_b :=\mathcal F$. Define the linear maps
\[
\begin{array}{lrll}
\phi: & \mathcal F_a \times \mathcal V_a \times \mathcal F_b \times  \mathcal V_b  & \to & \mathcal F_a \times \mathcal V_a \times \mathcal V_a \times \mathcal F_b \times \mathcal V_b \times \mathcal V_b \\[2mm]
& (f_a,v_a,f_b,v_b) & \mapsto & (f_a,v_a,v_a,f_b,v_b,v_b) \\[2mm]
\psi: & \mathcal F_a \times \mathcal V_a \times \mathcal F_b \times  \mathcal V_b  & \to & \mathcal F_a \times \mathcal F_b  \\[2mm]
& (f_a,v_a,f_b,v_b) & \mapsto & (f_a,f_b) \\[2mm]
\end{array}
\]
Then for the maximally monotone subspace
\[
\mathcal M_I = \{(v_a,v_b,v_a^*,v_b^*) \in \mathcal V_a \times \mathcal V_b \times \mathcal V^*_a \times \mathcal V^*_b \mid v_a=-v_b, v^*_a=v^*_b \},
\]
it can be readily checked that
\[
\mathcal M_a \circ \mathcal M_b = f_{\psi} (b_{\phi} \left(\mathcal M_a \times \mathcal M_I \times \mathcal M_b) \right),
\]
where  $\mathcal M_a \times \mathcal M_I \times \mathcal M_b$ is clearly maximally monotone. Then the proof finishes by applying Lemma~\ref{lem:pullback}.
\eproof

We immediately have the following corollary.
\begin{corollary}\label{cor:compo}
Consider a Dirac structure $\mathcal D \subset \mathcal X \times \mathcal F_R \times \mathcal X^* \times \mathcal F^*_R$, together with a nonnegative Lagrangian subspace $\mathcal R \subset \mathcal F_R \times \mathcal F^*_R$. Then the composition of $\mathcal D$ and $\mathcal R$ defined via
\[
\begin{array}{rcl}
\mathcal D \circ \mathcal R := \{(f,e) \in \mathcal X \times \mathcal X^* & \mid & \mbox{ there exists } f_R \in \mathcal F_R, e_R \in \mathcal F_R^*  \mbox{ such that  }\\[2mm]&&  (f,-f_R,e,e_R) \in \mathcal D, \mbox{ and }
 (f_R,e_R) \in \mathcal R \}
\end{array}
\]
is maximally monotone. In particular, for any $(f,e) \in \mathcal D \circ \mathcal R$, one has
\[
e^\top f = e_R^\top f_R \geq 0.
\]
%
\end{corollary}

\begin{remark}\label{rem:conj}{\rm
We conjecture that conversely any maximally monotone subspace $\mathcal M$ can be generated this way, i.e., as the composition of a certain Dirac structure $\mathcal D$ and a certain nonnegative Lagrangian subspace $\mathcal R$.
}
\end{remark}

The presented analysis of maximally monotone subspaces leads to the following geometric definition of an \emph{extended dHDAE} system, covering both dHDAE systems \eqref{gendhdae} and extended HDAE systems \eqref{gendila}.  See \cite{GerHR21} for related results (using the terminology of (maximally) dissipative linear relations).
%
\begin{definition}\label{def:edhdae}
Consider a linear state space $\mathcal X$ with coordinates $x$, a maximally monotone subspace $\mathcal M \subset \mathcal X \times \mathcal X^*$, and a Lagrange structure $\mathcal L \subset \mathcal X \times \mathcal X^*$. Then an \emph{extended dHDAE system}  is a system $(\mathcal X,\mathcal M,\mathcal L)$ satisfying
\begin{equation}\label{extdHDAE}
\{(\dot{x},x) \mid \mbox{ there exists } e \in \mathcal X^* \mbox{ such that } (-\dot{x},e) \in \mathcal M, (x,e) \in \mathcal L \}.
\end{equation}
\end{definition}
A coordinate representation of an extended dHDAE system is obtained as follows. Consider a coordinate expression \eqref{MN} of the maximally monotone subspace $\mathcal M$ (with $M,N$ satisfying \eqref{MN1}). This means that any element $(f,e) \in \mathcal M$ can be represented as
\[
\mat{c} f \\ e \rix = \mat{c} N^\top \\ M^\top \rix v
\]
for some $v \in \mathbb R^n$. Furthermore, any $(x,e) \in \mathcal L$ can be represented as in \eqref{SP1}. Substituting $-f=\dot{x}=P\dot{z}$ this yields
\[
\mat{c} -P\dot{z} \\ Sz \rix = \mat{c} N^\top \\ M^\top \rix v
\]
Now construct matrices $C,D$ satisfying
\bq
\label{CD}
\ker \mat{cc} C & D \rix = \Img \mat{c} N^\top \\ M^\top \rix.
\eq
Then pre-multiplication by such a maximal annihilator $\mat{cc} C & D \rix$ eliminates the auxiliary variables $v$, and one obtains the coordinate representation
\begin{equation} \label{CDham}
CP \dot{z} =DSz.
\end{equation}

\begin{remark}\label{rem:ports}{\rm
The geometric construction of extended dissipative Hamiltonian system can be immediately generalized to extended dissipative port-Hamiltonian DAE (dpHDAE) systems \emph{with external port variables} (inputs and outputs), by extending the maximally monotone subspace $\mathcal M \subset \mathcal X \times \mathcal X^*$ to a maximally monotone subspace $\mathcal M_e \subset \mathcal X \times \mathcal X^* \times \mathcal F_P \times \mathcal F_P^*$, where $\mathcal F_P \times \mathcal F_P^*$ is the space of external port variables.
}
\end{remark}

Two particular cases of Definition \ref{def:edhdae} are of special interest. The first one is where the maximally monotone subspace $\mathcal M$ is actually a Dirac structure as in \eqref{KL} with $K,L$ satisfying \eqref{dirac}. In this case one can take $\mat{cc} C & D \rix=\mat{cc} K & L \rix $, and thus the extended dHDAE system reduces to the extended HDAE system
\begin{equation}
\label{KLham}
KP \dot{z} = LSz.
\end{equation}
The other special case is where the maximally monotone subspace $\mathcal M$ in \eqref{MN} is such that $M$ is invertible. In this case, without loss of generality $M^\top$ can be taken to be the identity matrix, and the maximal annihilator $\mat{cc} C & D \rix$ can be taken to be of the form $\mat{cc} I & D \rix$. Hence $D=-N^\top$, and thus \eqref{CDham} reduces to
\[
P \dot{z} =-N^\top Sz.
\]
Furthermore, $MN^\top + NM^\top \geq 0$ reduces to $N^\top + N \geq 0$, and hence
\[
N^\top = - \frac{1}{2} (N - N^\top) + \frac{1}{2} (N + N^\top) =: -J + R
\]
with $J=-J^\top$ and $R=R^\top \geq 0$. Thus in this case the extended dpHDAE system takes the familiar form~\eqref{gendhdae} with $E=P$, $S=Q$ expressed as
\[
E\dot{z} = (J-R)Qz, \quad Q^\top E= E^\top Q, J=-J^\top, R=R^\top \geq 0.
\]
\begin{remark}\label{rem:constraints}{ \rm
Similar to the theory exposed in \cite{SchM18} for HDAE systems \eqref{gendila}, the \emph{algebraic constraints} of the dpHDAE system \eqref{CDham} can be split into two classes: one corresponding to singularity of $P$ (Lagrange algebraic constraints in \cite{SchM18}), and one corresponding to singularity of $C$. In case of \eqref{gendila} the second class of algebraic constraints are called Dirac algebraic constraints in \cite{SchM18}, but now they correspond to the maximally monotone subspace.

Furthermore, mimicking the developments in \cite{SchM18}, one can transform algebraic constraints associated with index one belonging to one class into algebraic constraints in the other, by the use of additional state variables (serving as Lagrange multipliers).
}
\end{remark}

\section{Representation of DAE systems generated by Dirac and Lagrange structures in the state variables $x$}\label{sec:geometric}

The representation \eqref{gendila} of an extended HDAE system as discussed in the previous sections does not use the state variables $x$ of the state space $\X$, but instead an equally dimensioned vector $z\in \mathcal Z$ parameterizing the Lagrange structure, cf. \eqref{imageL} and \eqref{lagrangespace}. In this section we show how a \emph{different} DAE representation involving the original state vector $x \in \X$ can be obtained. Furthermore we discuss in what sense this representation in $x$ is equivalent with the representation \eqref{gendila} involving $z$.

\subsection{A coordinate representation in the original state variables $x$}\label{sec:effortflow}
Consider a Dirac structure $\mathcal D \subset \mathcal X \times \X^*$, a Lagrange structure $\mathcal L \subset \mathcal X \times \mathcal X^*$, and the resulting dynamics specified (in coordinate-free form) as $\mathcal D \circ \mathcal L \subset {\mathcal X} \times \mathcal X$.
Let $x$ be coordinates for the state space $\mathcal X$ and let the Dirac structure represented by a pair of matrices $K,L$ and the Lagrange structure by a pair of matrices $P,S$.
To derive a coordinate representation employing directly the state vector $x$, we first consider the combined representations of $\mathcal D$ and $\mathcal L$, both in kernel representation, i.e.,
\begin{equation}
\label{eq:withe}
\mat{ccc} K & L & 0 \\ 0 & -P^\top & S^\top \rix \mat{c} -\dot{x} \\ e \\ x \rix=0,
\end{equation}
where $e$ are dual coordinates for $\X^*$. In order to obtain a DAE system only involving $x$ we need to eliminate the variables $e$. This can be done by considering a \emph{maximal annihilator} (left null-space) $\mat{cc} M & N \rix$ of $\mat{c} L \\ -P^\top \rix$, i.e.,
\begin{equation} \label{maxa}
\ker \mat{cc} M & N \rix = \Img \mat{c} L \\ -P^\top \rix,
\end{equation}
and thus, in particular,
\begin{equation}\label{MLNP}
ML = NP^\top .
\end{equation}
Since
\[
\mat{cc} M & N \rix \mat{ccc} K & L & 0 \\ 0 & -P^\top & S^\top \rix = \mat{ccc} MK & 0 & NS^\top \rix,
\]
premultiplication of the equations \eqref{eq:withe} by $\mat{cc} M & N \rix$ thus yields
\[
\mathcal D \circ \mathcal L = \ker \mat{cc} MK & NS^\top \rix\subset {\mathcal X} \times \mathcal X.
\]
Hence the resulting DAE system is given by
\begin{equation}\label{gendilax}
MK \dot{x} = NS^\top x.
\end{equation}
\begin{remark}{\rm
Also for extended dHDAE systems (including \emph{dissipation}) we can consider, instead of the coordinate representation \eqref{CDham} involving the parametrizing vector $z$, a representation that is using the original state $x$. In fact, let as before, cf. \eqref{CD}, $\mat{cc} C & D \rix$ denote a maximal annihilator of $\mat{cc} N & M \rix^\top$, i.e., $\ker \mat{cc} C & D \rix =\Img \mat{cc} N & M \rix^\top $. Then consider, similarly to \eqref{eq:withe}, the stacked matrix
\[
\mat{ccc} C & D & 0 \\ 0 & -P^\top & S^\top \rix
\]
and a maximal annihilator $\mat{cc} V & W \rix$ to $\mat{cc} D^\top & -P \rix^\top$, that is $\ker \mat{cc} V & W \rix = \Img \mat{cc} D^\top & -P \rix^\top$. Then premultiplication by $\mat{cc} V & W \rix$ yields the representation
\[
\label{gendilaxCD}
VC\dot{x}=WS^\top x
\]
The analysis performed in the current subsection for \eqref{gendilax} can be performed, mutatis mutandis, for \eqref{gendilaxCD} as well.
}
\end{remark}
%
%
%
%
%
%

Recall that in the coordinate representation \eqref{gendila} we have the expression $\cH^z(z)= \frac 12 z^\top S^\top Pz$ for the Hamiltonian. In the representation~\eqref{gendilax} we do not yet have  a Hamiltonian associated with the extended HDAE system. To
define such a Hamiltonian $\cH^x$ in \eqref{gendilax}  we would need that $P$ is invertible, in which case it is given by
\begin{equation} \label{defHx}
\cH^x(x)= \frac{1}{2}x^\top SP^{-1}x.
\end{equation}
If $P$ is invertible then there is a direct relation between the Hamiltonians~\eqref{defHzPS} and~\eqref{defHx}. In fact, substituting $x=Pz$, we immediately obtain
\[
\cH^x(x)= \frac{1}{2}x^\top SP^{-1}x= \frac{1}{2}z^\top P^\top SP^{-1}Pz = \frac{1}{2}z^\top P^\top Sz = \cH^z(z).
\]
Alternatively, if $S$ is invertible then one can
use the  \emph{co-energy} (Legendre transform) of $H^x$ given by
\begin{equation}\label{defHcoenergy}
\cH^e(e)= \frac{1}{2}e^\top PS^{-1}e
\end{equation}
for which $\cH^e(e) = \cH^z(z)$ with $e=Sz$.

Note that if $P$ is invertible, then also $M$ is invertible. This follows, since then $ML=NP^\top $ implies that the columns of $N$ are in $\Img M$, and since $\rank \mat{cc} M & N \rix =n$ this means $M$ is invertible. The converse that $M$ invertible implies $P$ invertible follows analogously. In a similar fashion, it follows  that $L$ is invertible if and only if $N$ is invertible.
These observations imply the following simplifications of the representation \eqref{gendilax} under additional assumptions.
\begin{itemize}
\item [1.]
\begin{itemize}
\item [a)] If $P$ is invertible  then $N=MLP^{-\top}$ and by multiplying \eqref{gendilax} from the left by $M^{-1}$ we obtain  the system
\begin{equation}\label{KLQ}
K \dot{x} = LP^{-\top}S^\top x = L \left(SP^{-1}\right)^\top x = L SP^{-1} x,
\end{equation}
where the last equality follows from $S^\top P = P^\top S$. This is exactly the form of a Hamiltonian DAE system in case of a general Dirac structure and a Lagrange structure that is given as the graph of a symmetric matrix $Q:=SP^{-1}$, see \cite{SchM95,Sch00,Sch13,SchJ14}. Indeed, the Lagrange structure simplifies to the gradient of the Hamiltonian function $\mathcal H^x(x)= \frac{1}{2}x^\top SP^{-1}x$.
\item [b)] If in this case additionally $K$ is invertible, then we obtain the \emph{Poisson formulation} of Hamiltonian systems, see e.g.   \cite{Arn13},
\begin{equation}\label{JQ}
\dot{x} = \left(K^{-1} L\right) \left(SP^{-1}\right) x= JQx,
\end{equation}
with $J=-J^\top =K^{-1} L$, and $Q=Q^\top $.
\end{itemize}
\item [2.]
\begin{itemize}
\item [a)] If  $L$ and thus also  $N$ is invertible, then by \eqref{MLNP} we have  $M=NP^\top L^{-1}$ and multiplying with $N^{-1}$ from the left we get the DAE
\begin{equation} \label{PJS}
P^\top J\dot x=P^\top L^{-1}K \dot{x} = S^\top x
\end{equation}
with $J:=\left(L^{-1}K \right)^\top = - L^{-1}K$.
\item [b)]
If additionally $P$ is invertible, then with $Q=SP^{-1}= P^{-\top} S^\top$ this may be rewritten as
\begin{equation}\label{JinvQ}
J\dot{x} = Qx,
\end{equation}
which is the standard \emph{symplectic formulation} of a Hamiltonian system in case additionally $J$ is invertible, see e.g.~\cite{Arn13}.
\end{itemize}
\end{itemize}

\subsection{Relation between the representations~\eqref{gendila} and~\eqref{gendilax}}

An immediate question that arises is how the representations~\eqref{gendila} and~\eqref{gendilax} are related. We have already seen that if $P$ is \emph{invertible} then the relationship is obvious, since in this case $x=Pz$ defines an ordinary state space transformation. However, if $P$ is not invertible then the representations are \emph{not} state space equivalent, as the following simple example demonstrates.

\begin{example}\label{ex:xnotz} {\rm For  $P=[0]$, $S=[1]$, $K=[1]$, $L=[0]$, we have that \eqref{gendila} is the singular system $0 \cdot \dot{z} = 0 \cdot z$. On the other hand
\[
\rank \mat{c} L \\ -P^\top \rix =0,
\]
and $\mathcal D \circ \mathcal L$ is the origin in ${\mathcal X} \times \mathcal X$, defining the degenerate DAE system $\dot{x}=x=0$.
}
\end{example}

However, representations~\eqref{gendila} (in the parameterizing $z$ variables) and~\eqref{gendilax} (in the original state variables $x$) can be shown to be equivalent in the following {\it generalized sense}.
First note that for any representation $P,S$ of a Lagrange structure there exist nonsingular matrices $V,W$ such that
\begin{eqnarray}
V^{-1}PW =V^\top  PW= \mat{cc} I_{n_1} & 0 \\ 0 & 0 \rix, \ V^\top SW = \mat{cc} S_{11} & 0 \\ 0 & I_{n_2} \rix,
\nonumber\\
KV=\mat{cc} K_{11} & K_{12} \\ K_{21} & K_{22} \rix,\ LV^{-T} =LV=\mat{cc} L_{11} & L_{12} \\ L_{21} & L_{22} \rix.
\label{redlagrane}
\end{eqnarray}
This is a direct consequence of Lemma~\ref{lem:PScanform} that will be presented in the next section.
Setting  $\mat{c} z_1 \\ z_2 \rix = W^{-1} z$ and $\mat{c} x_1 \\ x_2 \rix =V^{-1}x$, it follows that $z_1=x_1$ and $z_2=e_2$.
After such a transformation the system $KP \dot{z} = LSz$ takes the form
\begin{equation}
\label{eq:or}
\mat{cc} K_{11} & K_{12} \\ K_{21} & K_{22} \rix \mat{cc} I & 0 \\ 0 & 0 \rix \mat{c} \dot{x}_1 \\ \dot{e}_2 \rix =
\mat{cc} L_{11} & L_{12} \\ L_{21} & L_{22} \rix \mat{cc} S_{11} & 0 \\ 0 & I \rix \mat{c} x_1 \\ e_2 \rix.
\end{equation}
{If we \emph{add} to the vector $\mat{cc} x_1 \\ e_2 \rix$ the subvector $x_2$, and if we consider the equations \eqref{eq:or} {\it together}} with the original Lagrange algebraic constraint $x_2=0$, then the so extended system can be rewritten as
\begin{equation}
\label{eq:result}
\mat{cc} K_{11} & K_{12} \\ K_{21} & K_{22} \\0 & 0 \rix \mat{cc} I & 0  \\ 0 & 0  \rix \mat{c} \dot{x}_1 \\ \dot{x}_2 \rix =
\mat{cc} L_{11} S_{11} x_1 \\ L_{21} S_{11} x_1 \\ x_2 \rix +
\mat{c} L_{12} \\ L_{22} \\ 0 \rix e_2.
\end{equation}

On the other hand, as shown in Subsection~\ref{sec:effortflow}, the extended dHDAE system defined by the Dirac structure $\mathcal D$ and the Lagrangian structure $\mathcal L$ in the state space variables $x$ can be expressed as
\begin{equation}
\label{eq:result1}
\mat{cccccc} K_{11} & K_{12} & L_{11} & L_{12} & 0 & 0 \\
K_{21} & K_{22} & L_{21} & L_{22} & 0 & 0 \\
0 &0 & -I & 0 & S_{11} & 0 \\
0 &0 & 0 & 0 & 0 & I
\rix \mat{c} - \dot{x}_1 \\ - \dot{x}_2 \\ e_1 \\ e_2 \\ x_1 \\ x_2 \rix =0,
\end{equation}
with $e_1,e_2$ serving as \emph{auxiliary variables}. Instead of eliminating $e_1,e_2$ from these equations, as discussed in Subsection 4.1, we can only eliminate $e_1$ by premultiplication of \eqref{eq:result1} by the full row rank matrix
\[
\mat{cc} M & N \rix = \mat{cccc} I & 0 & L_{11} & 0 \\
0 & I  & L_{21} & 0 \\
0 & 0 & 0 & I
\rix,
\]
which directly leads to the system \eqref{eq:result}. This extended equivalence between \eqref{gendila} and \eqref{gendilax} is summarized in the following proposition.
\begin{proposition}\label{prop:equidila}
Consider the pHDAE representations \eqref{gendila} and \eqref{gendilax} defined by the same Lagrange structure $\mathcal L$ represented by matrices $P,S$, and by the same Dirac structure $\mathcal D$ represented by $K,L$. Consider a transformation such that $P,S$ and $K,L$ are transformed into the form \eqref{redlagrane} with corresponding partitioning
\[
x = \mat{c} x_1 \\ x_2 \rix, \; z = \mat{c} z_1 \\ z_2 \rix,
\]
where $x_1=z_1$. Adding to \eqref{gendila} the Lagrange algebraic constraint $x_2=0$  corresponding to $x=Pz$, the resulting dHDAE system is given by \eqref{eq:result}. This system is equivalent to the representation \eqref{eq:withe} of \eqref{gendilax} after elimination of the variables $e_1$.
\end{proposition}
Note that the subvector $e_2=z_2$ can be regarded as the \emph{Lagrange multiplier vector} corresponding to the constraint $x_2=0$. As such, $e_2=z_2$ does not contribute to the expression of the Hamiltonian $\mathcal H^z(z)$.

Let us illustrate the previous discussion with some further examples.

\begin{example} \label{exB} {\rm Consider a system in the form \eqref{gendila} with
$K= \mat{cc} I & 0 \\ 0 & I \rix$, $L= \mat{cc} 0 & I \\ -I & 0 \rix$, $P= \mat{cc} I & 0 \\ 0 & 0 \rix$, $S= \mat{cc} I & 0 \\ 0 & I \rix$ which is the general DAE
\[
\mat{cc} I & 0 \\ 0 & 0 \rix   \mat{c} \dot{z}_1 \\ \dot{z}_2  \rix= \mat{cc} 0 & I \\ -I & 0 \rix \mat{c} z_1 \\ z_2  \rix
\]
In order to compute the representation \eqref{gendilax}, we solve
\[
0= \mat{cc} M & N  \rix
\mat{cc}0 & I\\-I & 0 \\ -I & 0 \\ 0 & 0 \rix,
\]
and with
\[
\mat{cc} M & N  \rix = \mat{cccc} 0 & 0 & 0 & I  \\ 0 & I & -I & 0  \rix
\]
we get the system
\[
\mat{cc} 0 & 0 \\ 0 & I \rix   \mat{c} \dot{x}_1 \\ \dot{x}_2  \rix= \mat{cc} 0 & I \\ -I & 0 \rix \mat{c} x_1 \\ x_2  \rix.
\]
}
\end{example}

\begin{example}\label{exC}{\rm
Consider the system of the form \eqref{gendila} with $K= 0$, $L= I$, $P= I$, $S= I$, i.e.,
$0 \cdot \dot{z} = z$. Solving
\[
0= \mat{cc} M & N  \rix
\mat{c} I \\ -I \rix,
\]
yields $M=I, N=I$, and we obtain the system \eqref{gendilax} given by $0 \cdot \dot{x} = x$.
}
\end{example}

\begin{example}\label{exD1}{\rm
Consider the system of the form \eqref{gendila} with $K= 0$, $L= I$, $P= 0$, $S= I$, i.e. $0 \cdot \dot{z} = z$.
Solving
\[
0= \mat{cc} M & N  \rix
\mat{c}I \\ 0 \rix
\]
yields $M=0, N=I$, and thus the representation~\eqref{gendilax} is $0 \cdot \dot{x} = x$.
}
\end{example}

\begin{example}\label{exD}{\rm
Consider the system of the form \eqref{gendila} with $K= 0$, $L= \bma 0 & I \\ -I & 0 \ema$, $P= I$, $S= I$, i.e.,  $0 \cdot \dot{z} = \bma 0 & I \\ -I & 0 \ema z$.
Solving
\[
0= \mat{cc} M & N  \rix
\mat{cc}  0 & I \\ -I & 0 \\ -I & 0 \\ 0 & -I \rix
\]
yields $M= I, N= \mat{cc} 0 & I \\ -I & 0 \rix$, and thus the representation~\eqref{gendilax} is given by
\[
0 \cdot \dot{x} = \mat{cc} 0 & I \\ -I & 0\rix x.
\]
}
\end{example}

\section{Equivalence transformations and condensed forms}\label{sec:equi}
To characterize the properties  of extended dHDAEs we use transformations to condensed forms from which the properties can be read off.

For general DAEs~\eqref{gendae} given by matrix pairs $(E,A)$, $E,A\mathbb R^{\ell,n}$( or the representation via matrix pencils $\lambda E-A$) we can perform  equivalence transformations of the coefficients of the form
\begin{equation}\label{equisimple}
(\tilde E,\tilde A)=(U^\top EW,U^\top  AW),
\end{equation}
with $U\in \mathbb R^{\ell,ell}$, $W\in \mathbb R^{n,n}$ nonsingular. This corresponds to a scaling of the equation with $U^\top $ and a change of variables $z= W\tilde z$.
Under such transformations there is a one-to-one relationship between the solution spaces, see \cite{KunM06} and the canonical form is the Weierstra{\ss} canonical form.

For structured systems of the form \eqref{gendhdae}, the associated equivalence transformation that preserves the structure is of the form
\[
(\tilde E,(\tilde J-\tilde R)\tilde Q)=( U^\top EW,(U^\top  (J-R)U) U^{-1}QW),
\]
with $U\in \mathbb R^{\ell,\ell}$, $W\in \mathbb R^{n,n}$ nonsingular. A condensed form for this case has been presented in \cite{MehMW18}.

Finally for systems of the form \eqref{gendila}, the equivalence transformations have the form
\begin{equation}\label{equitrans}
\tilde K= U^\top KV,\ \tilde L= U^\top  L V^{-T},\ \tilde P = V^{-1} P W,\ \tilde S= V^\top  S W,
\end{equation}
where $U\in \mathbb R^{\ell,\ell}$, $V \in \mathbb R^{n,n}$, $W\in \mathbb R^{m,m}$ are nonsingular.

The geometric interpretation of the set of transformations in \eqref{equitrans} is clear: $V$ defines a coordinate transformation on the state space $\mathcal X$ while $V^{-\top}$ is the corresponding dual transformation on the dual state space $\mathcal X^*$. Also note that the combination of $V$ and $V^{-\top}$ on the product space $\mathcal X \times \mathcal X^*$ leaves the canonical bilinear forms defined by the matrices $\Pi_-$ and $\Pi_+$ invariant. (In fact, it can be shown that any transformation on $\mathcal X \times \mathcal X^*$ that leaves both canonical bilinear forms invariant is necessarily of this form for some invertible $V$.) Finally, $U^\top $ is an invertible transformation on the equation space for the kernel representation of the Dirac structure $\D$, while $W$ is an invertible transformation on the parametrization space $\mathcal Z$ for the Lagrange structure $\mathcal L$.

In all three cases, in view of an implementation of the transformations as numerically stable procedures, we are also interested in the case that $U,V,W$  are real orthogonal matrices. We then have that $V^{-1}=V^\top $ and for both pairs $(K,L)$ and $(P,S)$ this is a classical orthogonal equivalence transformation.

Using the described equivalence transformations we can derive condensed forms for pencils $\lambda P -S$ with $P^\top S= S^\top P$ associated with Lagrange subspaces (or \emph{isotropic subspaces}  if the dimension is not $n$, see e.g. \cite{MehMW18}). Here we slightly modify the representation and also give a constructive proof that can be implemented as numerically stable algorithm in Appendix A.
%
\begin{lemma}\label{lem:PScanform}
Let $P,S\in\mathbb R^{n,m}$ be such that $P^\top S=S^\top P$. Then there exist invertible matrices
$V\in\mathbb R^{n,n}$, $W\in\mathbb R^{m,m}$ such that
\begin{equation}\label{cformPS}
V^{-1}PW=\mat{ccccc}I_{m_1} & 0 & 0 & 0 & 0 \\ 0 & I_{m_2}&0 & 0 & 0 \\ 0 & 0 & I_{m_3} & 0 & 0 \\ 0 & 0 & 0 & 0 & 0 \\ 0 & 0 & 0 & 0 & 0 \\ 0 & 0 & 0 & 0 & 0\rix,
V^\top SW=\mat{cccccc}I_{m_1} & 0 & 0 & 0 & 0\\
0 & -I_{m_2} & 0 & 0 & 0\\
0 & 0 & 0 & 0 &  0\\
0 & 0 & 0 & I_{m_4} & 0 \\
S_{51} & S_{5,2} & S_{5,3} & 0  & 0\\
0 & 0 & 0 & 0 & 0
\rix,
\end{equation}
with $\mat{ccc} S_{51} & S_{5,2} & S_{5,3}\rix$ of full row rank $n_5$.
(Note that block sizes may be zero).
Moreover,  if the pencil $\lambda P-S$ is regular then the condensed form is unique, except for the order of blocks, and just contains the first four block rows and columns.
\end{lemma}
\proof
See Appendix A.
\eproof
Note that the condensed form is in general not unique in the fifth block row, but the block sizes $m_1,m_2,m_3,m_4$ and the row dimension $n_5$ are.
%
\begin{corollary}\label{cor:orthPS}
Let $P,S\in\mathbb R^{n,m}$ be such that $P^\top S=S^\top P$. Then there exist real orthogonal matrices
$V\in\mathbb R^{n,n}$, $W\in\mathbb R^{m,m}$ such that
\begin{equation}\label{condformPS}
V^\top PW=\mat{cccc}P_{11} & 0 & 0 & 0  \\ P_{21} & P_{22} &0 & 0  \\ 0 & 0  & 0 & 0 \\ 0 & 0 & 0 & 0  \\ 0 & 0 & 0 & 0  \rix,
V^\top SW=\mat{cccccc}S_{11} & 0 & 0 & 0 \\
0 & 0 & 0 & 0\\
S_{31} & S_{32} & S_{33} &  0 \\
S_{41} & S_{42} & 0  & 0
\\ 0 & 0 & 0 & 0
\rix,
\end{equation}
with $P_{11},S_{11}\in \mathbb R^{m_1+m_2,m_1+m_2}$, $P_{22} \in \mathbb R^{m_3,m_3}$, $S_{33} \in \mathbb R^{m_4,m_4}$ invertible, $\mat{ccc} S_{41} & S_{42} \rix$ of full row rank $n_5$, and $P_{11}^\top S_{11}=S_{11}^\top P_{11}$. Here the block sizes $m_1+m_2$, $m_3$, $m_4$, and $m_5$ are as in \eqref{cformPS}.
\end{corollary}
\proof
The proof follows by performing Steps 1. and 2. of the proof of Lemma~\ref{lem:PScanform}, see Appendix B, which yields
\[
V_2^\top  V_1^\top  P W_1 W_2= \mat{ccc} \hat P_{11} & 0 & 0\\ 0 & 0 & 0 \\ 0 & 0 & 0 \rix,\
V_2^\top  V_1^\top  S W_1 W_2= \mat{ccc} \hat S_{11} & 0 & 0\\ \hat S_{21} & \hat S_{22} & 0 \\ \hat S_{31} & 0 & 0 \rix,
\]
followed by a singular value decomposition  $\hat V^\top _3\hat S_{11} \hat W_3= \mat{cc} \check S_{11} & 0 \\ 0 & 0 \rix $
with $\check S_{11}$ nonsingular diagonal and a full rank decomposition $\check V_4^\top  \hat S_{31}\hat W_3 =\mat{cc}  S_{41} & S_{42} \rix$.
\eproof
Corollary~\ref{cor:orthPS} shows that the characteristic quantities $m_1+m_2$, $m_3$ and $m_4$, as well as $n_5$ can be obtained by purely real orthogonal transformations.   The quantities $m_1,m_2$ can then be determined from the real orthogonal staircase form of
the symmetric pencil $\lambda P_{11}^\top  S_{11} -S_{11}^\top  P_{11}$ which has been presented in \cite{ByeMX07} and implemented as production software in \cite{BruM08}.

There is an analogous condensed form for pencils of the form $\lambda K-L$ satisfying $KL^\top=-LK^\top$. For the case of regular pairs this directly follows from the canonical form presented in \cite{Cou90}, but again we present the construction so that it can be directly implemented as a numerical method, see Appendix B.
\begin{lemma}\label{lem:KLcanform}
Let $K,L\in\mathbb R^{\ell,n}$ be such that $KL^\top =-LK^\top $. Then there exist invertible matrices
$U\in\mathbb R^{\ell,\ell}$, $V\in\mathbb R^{n,n}$ such that
\begin{equation}\label{cformKL}
U^\top KV=\mat{cccccc}I_{\ell_1} & 0 & 0 & 0 & 0 & 0\\ 0 & I_{\ell_1}&0 & 0 & 0 & 0\\ 0 & 0 & I_{\ell_3} & 0 & 0 & 0\\ 0 & 0 & 0 & 0 & 0 & 0 \\
0 & 0 & 0 & 0  & 0 & 0\rix,
U^\top LV^{-\top}=\mat{cccccc} 0 & I_{\ell_1} & 0 & 0 & L_{15} & 0\\
 -I_{\ell_1} & 0 & 0 & 0 & L_{25} & 0\\
0 & 0 & 0 & 0 &  L_{35} & 0\\
0 & 0 & 0 & I_{\ell_4} & 0 & 0\\
0 & 0 & 0 & 0  & 0 & 0
\rix,
\end{equation}
with $\mat{c} L_{15} \\ L_{2,5} \\ L_{3,5}\rix$ of full column rank $n_5$. (Note that block sizes may be zero).
Moreover,  if the pencil $\lambda K-L$ is regular then the condensed form is unique except for the order of blocks and just contains the first four block rows and columns.
\end{lemma}
\proof
See Appendix B.
\eproof

Note again that the form \eqref{cformKL} is not unique in general but the block sizes  $\ell_1,\ell_2,\ell_3,\ell_4$ and the column  dimension $n_5$ are.

\begin{corollary}\label{cor:orthKL}
Let $K,L\in\mathbb R^{\ell,n}$ be such that $LK^\top =-KL^\top $. Then there exist real orthogonal matrices
$U\in\mathbb R^{\ell,\ell}$, $V\in\mathbb R^{n,n}$ such that
\begin{equation}\label{condformKL}
U^\top KV=\mat{ccccc}K_{11} & K_{12} & 0 & 0 &0 \\ 0 & K_{22} &0 & 0  &0\\ 0 & 0  & 0 & 0 &0\\ 0 & 0 & 0 & 0  &0 \rix,
U^\top LV=\mat{cccccc}L_{11} & 0 & L_{13} & L_{14} &0  \\
0 & 0 &  L_{23} & L_{24} & 0\\
0 & 0 & L_{33} &  0  & 0\\
 0 & 0 & 0 & 0 & 0
\rix.
\end{equation}
with $K_{11},L_{11}\in \mathbb R^{2\ell_1,2\ell_1}$, $K_{22} \in \mathbb R^{\ell_3,\ell_3}$, $L_{33} \in \mathbb R^{\ell_4,\ell_4}$ invertible, $\mat{c} L_{14} \\ L_{2,4} \rix$ is of full column rank $n_5$, and $K_{11}^\top L_{11}=-L_{11}^\top K_{11}$. Here the block-sizes $\ell_1$, $\ell_3$, $\ell_4$, and $n_5$ are as in \eqref{cformKL}.
\end{corollary}
\proof
The proof follows by performing Steps 1. and 2. of the proof of Lemma~\ref{lem:KLcanform}, which yields
\[
U_2^\top  U_1^\top  K V_1 V_2= \mat{ccc} \hat K_{11} & 0 & 0\\ 0 & 0 & 0 \\ 0 & 0 & 0 \rix,\
U_2^\top  U_1^\top  L V_1 V_2= \mat{ccc} \hat L_{11} & \hat L_{12} & \hat L_{13}\\ 0 & \hat L_{22} & 0\\ 0 & 0 & 0 \rix,
\]
followed by a singular value decomposition  $\hat U^\top _3\hat L_{11} \hat V_4= \mat{cc} \check L_{11} & 0 \\ 0 & 0 \rix $
with $\check L_{11}$ nonsingular diagonal and a full rank decomposition $\hat U_3^\top \hat L_{13}\hat V_3 =\mat{c}  L_{14} \\ L_{24} \rix$.
\eproof
Corollary~\ref{cor:orthKL} shows that the characteristic quantities $\ell_1$, $\ell_3$, $\ell_4$, as well as $n_5$ can be obtained by purely real orthogonal transformations.

The presented condensed forms  can now be used in generating a condensed form for systems of the form \eqref{gendila}.
%
\begin{lemma}\label{lem:step1}
Consider a system of the form \eqref{gendila} with $K,P,L,S\in \mathbb R^{n,n}$ and regular pencil $\lambda KP-LS$. Then there exists invertible matrices $U,V,W$ as in \eqref{equitrans} such that
\begin{eqnarray}
\hat K&=& U^\top KV =\mat{cc} K_{11} & K_{12} \\ K_{21} & K_{22} \rix,\  \hat L =U^\top LV^{-\top} = \mat{cc} L_{11} & 0 \\ 0 & I_{n_2} \rix,\nonumber \\
 \hat P&=&V^{-1} P W=\mat{cc} I_{n_1} & 0 \\ 0 & 0 \rix,\ \hat S = V^\top  S W=\mat{cc} S_{11} & 0 \\ 0 & I_{n_2} \rix,
 \label{firstsimplerform}
\end{eqnarray}
where $S_{11}=S_{11}^\top $ and  $K_{11} L_{11}^\top =-L_{11} K_{11}^\top $.
\end{lemma}
\proof
Since the pencil $\lambda KP-LS$ is square and regular, it is square, and also the pencil $\lambda P-S$ is regular, otherwise by Lemma~\ref{lem:PScanform} there would be common kernel of $P$ and $S$ which would imply the pencil $\lambda KP-LS$ to be singular.

Thus, by Lemma~\ref{lem:PScanform} there exist nonsingular matrices
$W_1, V_1\in\mathbb R^{n,n}$  such that
\begin{eqnarray*}
\tilde K&=&KV_1=\mat{cc} \tilde K_{11} & \tilde K_{12} \\ \tilde K_{21} & \tilde K_{22} \rix,\
\tilde L=L V_1^{-\top}=\mat{cc} \tilde L_{11} & \tilde L_{12} \\ \tilde L_{21} & \tilde L_{22} \rix,\\
\tilde P&=&V_1^{-1} P W_1= \mat{cc} I_{n_1} & 0 \\ 0 & 0 \rix,\
\tilde S= V_1^\top  S W_1 =\mat{cc} \tilde S_{11} & 0 \\ 0 & I_{n_2}\rix,
\end{eqnarray*}
with $n_1= m_1+m_2+m_3$, $n_2 =m_4$ and $\tilde S_{11}$ symmetric. The regularity of the pencil $\lambda KP-LS$ implies that
\[
\mat{c} \tilde L_{12} \\ \tilde L_{22} \rix
\]
has full column rank and hence there exist invertible matrices $U_2\in \mathbb R^{n,n}$, $\tilde V_2\in \mathbb R^{n_2,n_2}$,  and
\[
V_2=\mat{cc} I_{n_1} & 0 \\ 0 & \tilde V_2\rix
\]
such that
\[
U_2^\top  \tilde L V_2=\mat{cc} \hat L_{11} & 0 \\ \hat L_{21} & I_{n_2} \rix.
\]
With
\[
V_3=\mat{cc} I_{n_1} & \hat L_{21}^\top  \\ 0 & I_{n_2} \rix,
\]
we then get that
\begin{eqnarray*}
\hat K&=& U_2^\top KV_1V_2V_3 =\mat{cc} K_{11} & K_{12} \\ K_{21} & K_{22} \rix,\  \hat L =U_2^\top LV_1^{-\top}V_2^{-\top}V_3^{-\top} = \mat{cc} L_{11} & 0 \\ 0 & I_{n_2} \rix,\\
\hat P&=&V_3^{-1}V_2^{-1}V_1^{-1} P W_1= \mat{cc} I_{n_1} & 0 \\ 0 & I_{n_2} \rix,\ \hat S = V_3^\top  V_2^\top V_1^\top  SW_1 =\mat{cc} S_{11} & 0 \\ 0 & I_{n_2} \rix,
\end{eqnarray*}
has the desired form with $U=U_2$, $V=V_1V_2V_3$, $W=W_1$, and where $S_{11}=S_{11}^\top $ and  $K_{11} L_{11}^\top =-L_{11} K_{11}^\top $.
\eproof

Transforming the system as in \eqref{firstsimplerform} and setting
\[
 z= W \mat{c} z_1 \\ z_2 \rix,
\]
partitioned accordingly, from the first block row of the coefficient matrices we obtain a reduced system given by
\begin{equation}\label{redP=I}
\bar K \bar P \dot z_1 = \bar L \bar S z_1,
\end{equation}
with $\bar P=I_{n_1}$, $\bar S =\bar S^\top  =S_{11}$, $\bar K=K_{11}$ and $\bar L= L_{11}$, together with an equation $z_2=K_{21} \dot z_1$, where $z_2$ does not contribute to the Hamiltonian $\mathcal H^z(z)=\frac 12 z^\top P^\top S z$. Note that the second equation  is an index two constraint, because it uses the derivative of $z_1$, \cite{KunM06}. It arises from the Lagrange structure due to the singularity of $P$.

An analogous representation can be constructed from the  condensed form of Lemma~\ref{lem:KLcanform}.
\begin{lemma}\label{lem:ste21}
Consider a system of the form \eqref{gendila} with $K,P,L,S\in \mathbb R^{n,n}$ and regular pencil $\lambda KP-LS$. Then there exist invertible matrices $U,V,W$ as in \eqref{equitrans} such that
\begin{eqnarray}
\hat K&=& U^\top KV =\mat{cc} I_{n_1} & 0 \\ 0 & 0 \rix,\  \hat L =U^\top LV^{-\top} = \mat{cc} L_{11} & 0 \\ 0 & I_{n_2} \rix,\nonumber \\
 \hat P&=&V^{-1} P W=\mat{cc} P_{11} & P_{12} \\ P_{21} & P_{22} \rix,\ \hat S = V^\top  S W=\mat{cc} S_{11} & 0 \\ 0 & I_{n_2} \rix,
 \label{secondsimplerform}
\end{eqnarray}
where $L_{11}=-L_{11}^\top $ and  $P_{11}^\top  S_{11}=S_{11}^\top  P_{11}$.
\end{lemma}
\proof
Since the pencil $\lambda KP-LS$ is square and regular, also the pencil $\lambda K-L$ is regular, otherwise by Lemma~\ref{lem:KLcanform} there would be a common left nullspace of $K$ and $L$ which would imply the pencil $\lambda KP-LS$ to be singular.

Thus by Lemma~\ref{lem:KLcanform}, there exist nonsingular matrices
$U_1, V_1\in\mathbb R^{n,n}$  such that
\begin{eqnarray*}
\tilde K&=&U_1KV_1=\mat{cc} I_{n_1} & 0 \\ 0 & 0 \rix,\
\tilde L=U_1 L V_1^{-\top}=\mat{cc} \tilde L_{11} & 0 \\ 0 & I_{n_2} \rix,\\
\tilde P&=&V_1^{-1} P = \mat{cc} \tilde P_{11} & \tilde P_{12} \\ \tilde P_{21} & \tilde P_{22} \rix,\
\tilde S= V_1^\top  S  =\mat{cc} \tilde S_{11} & \tilde S_{12} \\ \tilde S_{21} & \tilde S_{22}\rix,
\end{eqnarray*}
with $n_1= 2\ell_1+\ell_3$, $n_2 =\ell_4$ and $\tilde L_{11}$ skew-symmetric. The regularity of the pencil $\lambda KP-LS$ implies that
\[
\mat{cc} \tilde S_{21} &\tilde S_{22} \rix
\]
has full row rank and hence there exist invertible matrices $W_2\in \mathbb R^{n,n}$, $\tilde V_2\in \mathbb R^{n_2,n_2}$,  and
\[
V_2=\mat{cc} I_{n_1} & 0 \\ 0 & \tilde V_2\rix
\]
such that
\[
V_2^\top  \tilde S W_2=\mat{cc} \hat S_{11} & \hat S_{12} \\ 0 & I_{n_2} \rix.
\]
With
\[
V_3=\mat{cc} I_{n_1} & - \hat S_{12} \\ 0 & I_{n_2} \rix,
\]
we then get that
\begin{eqnarray*}
\hat K&=& U_1^\top KV_1V_2V_3 =\mat{cc} I_{n_1} & 0 \\ 0 & 0 \rix,\  \hat L =U_1^\top L V_1^{-\top}V_2^{-\top}V_3^{-\top} = \mat{cc} L_{11} & 0 \\ 0 & I_{n_2} \rix,\\
\hat P&=&V_3^{-1}V_2^{-1}V_1^{-1} P W_2= \mat{cc} P_{11} & P_{12} \\ P_{21} & P_{22} \rix,\ \hat S = V_3^\top  V_2^\top V_1^\top  S W_2 =\mat{cc} S_{11} & 0 \\ 0 & I_{n_2} \rix,
\end{eqnarray*}
has the desired form with $W=W_2$, $V=V_1V_2V_3$, $U=U_1$, and where $L_{11}=-L_{11}^\top $ and  $P_{11}^\top  S_{11}=S_{11}^\top  P_{11}$.
\eproof

Transforming $KP\dot z=LSz$ as in \eqref{secondsimplerform} and setting
\[
 z= W \mat{c} z_1 \\ z_2 \rix,
\]
partitioned accordingly, from the first block row of the coefficient matrices we obtain a reduced system given by
\begin{equation}\label{redK=I}
\bar K \bar P \dot z_1 = \bar L \bar S z_1,
\end{equation}
with $\bar P=P_{11}$, $\bar S =S_{11}$, $\bar K=I_{n_1}$ and $\bar L= L_{11}=-L_{11}^\top$, together with a differential algebraic equation $0=z_2$, so that $z_2$ does not contribute to the Hamiltonian.

\begin{remark}\label{rem:reduced}{\rm
System~\eqref{redP=I} is a dHDAE  of the form \eqref{KLQ} in which the Lagrange structure is spanned by the columns of
\[
\mat{c} I_{n_1} \\ \bar S_{11} \rix.
\]
See also \cite{MehMW21,MehU23} for similar constructions in the context of removing the factor $Q$ in systems of the form \eqref{gendhdae}.

Similarly, System~\eqref{redK=I} is a pHDAE of the form \eqref{gendhdae} (with $R=0$) and the Dirac structure is spanned by the columns of
\[
\mat{c} I_{n_1} \\ \bar L_{11} \rix.
\]
%
}
\end{remark}


We also perform a similar construction for systems  of the form \eqref{CDham}. Since $C$ and $D$ are chosen to be a maximal annihilator such that $C N^\top + D M^\top=0$ in \eqref{CD} and $MN^\top+NM^\top\geq 0$  with $\rank \mat{cc} N & M \rix =n$ we can use the same construction as in the proof of Lemma~\ref{lem:KLcanform} to first transform $N$ and $M$ in such a way that
\[
N=\mat{cc} I_{n_1} & 0 \\ 0 & 0 \rix,\ M=\mat{cc} M_{11} & 0 \\ 0 & I_{n_2} \rix.
\]
This implies that we may choose $C$ and $D$ such that
\[
C=\mat{cc} -M_{11}^\top & 0 \\ \\ 0 & I_{n_2} \rix, \ D=\mat{cc} I_{n_1} & 0 \\ 0 & 0 \rix.
\]
If $\lambda CP-DS$ is regular, then it follows that the last $n_2$ rows of $CP$ have full row rank and hence
altogether we have the following condensed form.
\begin{lemma}\label{lem:step1cd}
Consider a system of the form \eqref{CDham} with $C,P,L,S\in \mathbb R^{n,n}$ and regular pencil $\lambda CP-DS$ and $\mat{cc} C& D\rix$ a maximal annihilator as in \eqref{CD}. Then there exist invertible matrices $U,V,W$ as in \eqref{equitrans} such that
\begin{eqnarray}
\hat C&=& U^\top CV =\mat{cc} C_{11} & 0 \\ 0 & I_{n_2} \rix,\  \hat D =U^\top DV^{-\top} = \mat{cc} I_{n_1} & 0 \\ 0 & 0\rix,\nonumber \\
 \hat P&=&V^{-1} P W=\mat{cc} P_{11} & P_{12} \\ 0 & I_{n_2} \rix,\ \hat S = V^\top  S W=\mat{cc} S_{11} & S_{12} \\ S_{21} & S_{22} \rix,
 \label{cdsimplerform}
\end{eqnarray}
where $C_{11}=-C_{11}^\top $ and $P_{11}^\top S_{11}= S_{11}^\top P_{11}$.
\end{lemma}

As a consequence, by transforming $CP\dot z=DSz$ as in \eqref{cdsimplerform} and setting
\[
 z= W \mat{c} z_1 \\ z_2 \rix,
\]
partitioned accordingly, from the second block row of the coefficient matrices we obtain
$\dot z_2=0$, i.e. $z_2$ is a constant function and the first block row  gives an inhomogeneous reduced system
\begin{equation}\label{redD=I}
\bar C \bar P \dot z_1 = \bar D \bar S z_1 + S_{12} z_2,
\end{equation}
with $\bar P=P_{11}$, $\bar S =S_{11}$, $\bar P^\top\bar S=\bar S^\top \bar P$, $\bar C=-M_{11}^T$ and $\bar D= I_{n_1}$. It will then depend on the initial condition for $z_2$ whether $z_2=0$ in which case it does not contribute to the Hamiltonian, otherwise the Hamiltonian is still a quadratic function in $z_1$ plus some linear and constant terms.

\begin{remark}\label{rem:index2}{\rm
The condensed forms in this section require rank decisions. Even if they are done in a numerically stable way using singular value decompositions, they can give wrong decisions in finite precision arithmetic. It is a common strategy to use in the case of doubt the worst case scenario. In the case of condensed forms this would be to assume that the problem is a DAE of index two.
}
\end{remark}

In this section we have derived structured condensed forms and shown that these can also be used to identify a subsystem which is of one of the well-established forms plus an algebraic constraint whose solution does not contribute to the Hamiltonian. In the next section we analyze, when general DAEs can be transformed to the forms \eqref{gendhdae} or \eqref{gendila}.

\section{Representation of DAEs into the form $KP\dot z=LS z$ or $\frac{d}{dt}(Ez) =(J-R)Qz$}\label{sec:daetoph}
For general DAE systems $E\dot x=Ax$ it has been characterized in \cite{MehMW21} when they  are equivalent to a dHDAE system of the form \eqref{gendhdae}. We present here a simplified  result for the regular case.
\begin{theorem}\label{stable=dh}
~\\
i)  A regular pencil $L(\lambda)=\lambda \hat E-\hat A$ is equivalent to a pencil of the form
$\lambda E-(J-R)Q$ as in~\eqref{gendhdae} with $\lambda E-Q$ being regular
if and only if the following conditions are satisfied:
\begin{enumerate}
\item The spectrum of $L(\lambda)$ is contained in the closed left half plane.
\item The finite nonzero eigenvalues on the imaginary axis are semisimple and the partial multiplicities of the eigenvalue zero are at most two.
\item The index of $L(\lambda)$ is at most two.
\end{enumerate}

ii) A regular pencil $L(\lambda)=\lambda \hat E-\hat A$ is  equivalent to a pencil of
the form $\lambda E-(J-R)$ as in~\eqref{gendhdae} (i.e., with $Q=I$) if and only if the following conditions are satisfied:
\begin{enumerate}
\item The spectrum of $L(\lambda)$ is contained in the closed left half plane.
\item The finite eigenvalues on the imaginary axis (including zero) are semisimple.
\item The index of $L(\lambda)$ is at most two.
\end{enumerate}
\end{theorem}

As a Corollary for the case without dissipation we have the following result.

\begin{corollary}\label{cor:lossless}
A regular pencil $L(\lambda) =\lambda \hat E -\hat A$ is  equivalent to a pencil of
the form $\lambda E-J$ as in~\eqref{gendhdae} (with $Q=I,R=0$) if and only if the following conditions are satisfied:
\begin{enumerate}
\item All finite eigenvalues are on the imaginary axis  and semisimple.
\item The index of $L(\lambda)$ is at most two.
\end{enumerate}
\end{corollary}

To study when general regular DAEs of the form~\eqref{gendae} can be expressed as extended dHDAEs of the form \eqref{gendila} we first consider a condensed form under orthogonal equivalence.

\begin{theorem}\label{thm:sinind2}
Consider a regular pencil $\lambda E-A$ with $E,A\in \mathbb R^{n,n}$ of index at most two. Then there exist real orthogonal matrices $U\in \mathbb R^{n,n}$ and $V\in \mathbb R^{n,n}$ such that
\begin{equation}\label{cformind2}
U^\top  EV= \mat{cccc}  E_{11} & E_{12} & 0 & 0  \\
E_{21} & E_{22} & 0 & 0  \\
0& 0 & 0 & 0  \\
0 & 0 & 0 & 0
\rix,\ U^\top  A V=\mat{ccccc}  A_{11} & A_{12} & A_{13} & A_{14}   \\
A_{21} & A_{22} & A_{23}  & 0  \\
A_{31} & A_{32} & A_{33}  & 0  \\
A_{41} & 0 & 0 & 0
\rix,
\end{equation}
with $A_{14}\in \mathbb R^{n_1,n_1}$, $A_{41}\in \mathbb R^{n_1,n_1}$,
$E_{22}\in \mathbb R^{n_2,n_2}$, and $A_{33}\in \mathbb R^{n_3,n_3}$ invertible.
\end{theorem}
\proof
The proof is presented in Appendix C.
\eproof

Transforming the DAE \eqref{gendae} as  $U^\top EV V^\top \dot x= U^\top  A V V^\top x$ and setting
$V^\top  x= [x_1^\top , \ldots, x_4^\top ]^\top $, it follows that $x_1=0$, $x_2$ is determined form the  implicit ordinary differential equation (note that $E_{22}$ is invertible)
\begin{equation}\label{redDAE}
E_{22} \dot x_2 = (A_{22}-A_{23} A_{33}^{-1} A_{32}) x_2,
\end{equation}
$x_3= -A_{33}^{-1} A_{32} x_2$, and  $x_4$ is uniquely determined in terms of  $x_2,\dot x_2, x_3$.
Initial conditions can be prescribed freely for $x_2$ only.
\begin{corollary}\label{cor:canforma}
Consider a general regular pencil $\lambda E-A$  with $E,A\in \mathbb R^{n,n}$ that is of index at most two and for which all finite eigenvalues are in the closed left half plane and those on the imaginary axis are semi-simple. Then  there exist invertible matrices $U\in \mathbb R^{n,n}$ and $V\in \mathbb R^{n,n}$ such that
\begin{equation}\label{cformindred}
U^\top  EV= \mat{ccccc}  I_{\hat n_1} & 0 & 0  & 0  & 0\\
0 & 0 & 0 & 0 & 0\\
0 & 0 & E_{33} & 0& 0  \\
0 & 0 & 0 & 0  & 0\\
0 & 0 & 0 & 0  & 0\\
0 & 0 & 0 & 0 & 0
\rix,\ U^\top  A V=\mat{ccccc} 0& 0 & 0 &  I_{\hat n_1} &0 \\
 0 & 0 & 0 & 0 & I_{\hat n_2} \\
0 & 0 & A_{33} & 0  & 0 \\
0 & 0 &  0 &I_{\hat n_4}  & 0  \\
-I_{\hat n_1} & 0 & 0 & 0  & 0 \\
0 & -I_{\hat n_2} & 0 & 0 & 0
\rix,
\end{equation}
where
\[
E_{33}=E_{33}^\top>0,\ A_{33}= J_{33}-R_{33},\  J_{33}=-J_{33}^\top,\ R_{33}^\top=R_{33}\geq 0.
\]
\end{corollary}
\proof
The proof follows by considering the condensed form~\eqref{cformind2}, and using block elimination with the invertible matrices $A_{33}$, $A_{14}$, $A_{41}$, $E_{22}$ to transform pencil in~\eqref{cformind2} to
the form
\begin{equation}\label{almostph}
\lambda \mat{cccc}  \tilde E_{11} & 0 & 0 & 0   \\
0 & I_{n_2} & 0& 0   \\
0 & 0 & 0 & 0   \\
0 & 0 & 0 & 0
\rix-\mat{ccccc} 0& 0 & 0 &  I_{n_1}  \\
0 & \tilde A_{22} & 0  & 0 \\
0 & 0 & I_{n_3}  & 0  \\
-I_{n_1} & 0 & 0 & 0
\rix.
\end{equation}
Let  $U_1^\top \tilde E_{11} V_1= \mat{cc} I_{\hat{n}_1}& 0 \\ 0 & 0 \rix$ be the echelon form of $\tilde E_{11}$. We scale the first block row of with $U_1^\top$, the fourth block row with $V_1^{-1}$, the first block column by $V_1$, and the fourth block column by $U_1^{-\top}$ and obtain a form
\[
U^\top  EV= \mat{cccccc}  I_{{\hat n}_1} & 0 & 0  & 0  & 0 &0\\
0 & 0 & 0 & 0 & 0 & 0\\
0 & 0 & I_{n_3} & 0& 0 & 0 \\
0 & 0 & 0 & 0  & 0& 0\\
0 & 0 & 0 & 0  & 0& 0\\
0 & 0 & 0 & 0  & 0& 0
\rix,\ U^\top  A V=\mat{cccccc} 0& 0 & 0 & 0 & I_{{\hat n}_1} &0 \\
 0 & 0 & 0 & 0 &  0&  I_{{\hat n}_2} \\
0 & 0 & \hat A_{33} & 0  & 0 & 0\\
0 & 0 &  0 &I_{{\hat n}_4}  & 0 & 0 \\
-I_{{\hat n}_1} & 0 & 0 & 0  & 0 & 0 \\
0 & -I_{{\hat n}_2} & 0 & 0 & 0 & 0
\rix,
\]

For any positive definite solution $X$ of the Lyapunov inequality
\begin{equation}\label{lyaineq}
-\tilde A_{33}^\top  X-X \tilde A_{33}\geq 0
\end{equation}
one can multiply the second block row by $X$ and obtain  that $E_{33}=X$ and $A_{33} =X\hat A_{33}$ has the desired form, see e.g. \cite{AchAM21,BeaMV19}.
\eproof
Note that the transformation to a system of the form \eqref{gendhdae} can also be achieved in a similar way for singular pencils with zero minimal indices.

If there is no dissipation, i.e. if $R_{33}=0$, then $A_{33}$ is skew-symmetric.

In Corollary~\ref{cor:canforma} we have shown that general systems $E\dot x=Ax$  can be transformed
to a very special canonical form and the following remark shows that for the case $R_{33}=0$ each of the blocks in the canonical form can be expressed as a pencil of the form $\lambda KP-LS$.

\begin{remark}\label{rem:klpsrepresentation}{\rm
Consider a regular pencil $\lambda E-A$ in the form \eqref{cformindred} then after a permutation one gets four blocks which all can be written in the form $\lambda KP-LS$ as in \eqref{gendila}.
\begin{itemize}
\item [1)] We have
\[ \mat{cc} I_{{\hat n}_1} & 0 \\ 0 & 0 \rix \mat{c} \dot z_1 \\ \dot z_5 \rix- \mat{cc} 0 & I_{{\hat n}_1} \\ -I_{{\hat n}_1} & 0 \rix\mat{c} z_1 \\ z_5 \rix=
 K_{1}P_{1}\mat{c} \dot z_1 \\ \dot z_5 \rix- L_{1}S_{1}\mat{c} z_1 \\ z_5 \rix
\]
 with
 \[
 K_1=\mat{cc} I_{{\hat n}_1} & 0 \\ 0 & I_{{\hat n}_1} \rix,\ L_1= \mat{cc} 0 & I_{{\hat n}_1} \\ -I_{{\hat n}_1} & 0 \rix,\
P_1= \mat{cc} I_{{\hat n}_1} & 0 \\ 0 & 0 \rix, \ S_1=  \mat{cc} I_{{\hat n}_1} & 0 \\ 0 & I_{{\hat n}_1}  \rix
\]
and Hamiltonian $\mathcal H^z_1= \frac 12 \mat{c} z_1 \\ z_5\rix^\top  \mat{cc} I_{{\hat n}_1} & 0 \\ 0 & 0 \rix \mat{c} z_1 \\ z_5 \rix={\frac 12 z_1^\top z_1}$ which is actually $0$.

In this case  we can  insert the derivative of the second equation into the first (index reduction) and obtain
%
%
%
\[
\mat{cc} 0 & 0 \\ 0 & 0 \rix \mat{c} \dot z_1 \\ \dot z_5 \rix= \mat{cc} 0 & I_{{\hat n}_1} \\ -I_{{\hat n}_1} & 0 \rix \mat{c} z_1 \\ z_5 \rix
\]
without changing the Hamiltonian.

\item[2)] We have
\[
\mat{cc} 0 & 0 \\ 0 & 0 \rix \mat{c} \dot z_2 \\ \dot z_6 \rix=\mat{cc} 0 & I_{{\hat n}_2} \\ -I_{{\hat n}_2} & 0 \rix\mat{c} z_2 \\ z_6 \rix=
K_{2}P_{2}\mat{c} \dot z_2 \\ \dot z_6 \rix- L_{2}S_{2}\mat{c} z_2 \\ z_6 \rix
\]
with different possibilities of representation, e.g.\\
a)  \[
 K_2=\mat{cc} I_{{\hat n}_2} & 0 \\ 0 & I_{{\hat n}_2} \rix,\ L_3= \mat{cc} 0 & I_{{\hat n}_2} \\ -I_{{\hat n}_2} & 0 \rix,\
P_2= \mat{cc} 0 & 0 \\ 0 & 0 \rix, \ S_3=  \mat{cc} I_{{\hat n}_1} & 0 \\ 0 & I_{{\hat n}_1}  \rix
\]
and Hamiltonian $\mathcal H^z_2= \frac 12 \mat{c} z_2 \\ z_6\rix^\top  \mat{cc} 0 & 0 \\ 0 & 0 \rix \mat{c} z_2 \\ z_6 \rix=0$, or

b)  \[
 K_2=\mat{cc} 0 & 0 \\ 0 & 0 \rix,\ L_2= \mat{cc} 0 & I_{{\hat n}_2} \\ -I_{{\hat n}_2} & 0 \rix,\
P_2= \mat{cc} 0 & 0 \\ 0 & 0 \rix, \ S_2=  \mat{cc} I_{{\hat n}_2} & 0 \\ 0 & I_{{\hat n}_2}  \rix
\]
and Hamiltonian $\mathcal H^z_2= \frac 12 \mat{c} z_2 \\ z_6\rix^\top  \mat{cc} 0 & 0 \\ 0 & 0 \rix \mat{c} z_2 \\ z_6 \rix=0$.
\item [3)] We have $\lambda E_{33}-A_{33}= \lambda K_{3}P_{3}- L_{3}S_{3}$ with
 $K_3=I_{\hat{n_3}}$, $P_3=E_{33}$, $L_3= A_{33}$, $S_{3}= I_{{\hat n}_3}$.
Here the Hamiltonian is $\mathcal H^z_3= \frac 12 z_3^\top  E_{33} z_3$.
\item [4)] We have $\lambda 0 -I_{{\hat n}_4}= \lambda K_{4}P_{4}- L_{4}S_{4}$
with $K_4=0$, $P_4=I_{{\hat n}_4}$, $L_4=I_{{\hat n}_4}$, $S_4=I_{{\hat n}_4}$. Here the Hamiltonian is $\mathcal H^z_4= \frac 12 z_4^\top  z_4$.
\end{itemize}

Note that the presented representations are in no way unique, but if the condensed form is available or computable, and the properties of Theorem~\ref{thm:sinind2} hold then we can express the general DAE in the representation \eqref{gendila} or \eqref{gendhdae}.
}\end{remark}
{ This discussion yields the following useful corollary.
\begin{corollary}\label{cor:lagrangeind2}
Consider a regular pencil of the form $\lambda KP-LS$ associated with the dHDAE \eqref{gendila}. Then it  has index at most two, and index two can only occur if the system has a singular Lagrange structure.
\end{corollary}
\proof
Consider the representations in Remark~\ref{rem:klpsrepresentation}. Then the index two structure occurs only in the first case where $K_1,L_1$ are invertible, but the product $P_1^TS_1$ is singular. Thus index two arises only from a singular Lagrange structure.
\eproof
}
\section*{Conclusion and Outlook}

Different definitions of (extended, dissipative) Hamiltonian or port Hamiltonian differential-algebraic systems lead to different representations. We have collected all the known representations as well as a few new ones and analyzed them from a geometric as well as an algebraic point of view. The latter leads to condensed forms  that can be directly implemented in numerical algorithms to compute the structural properties of the systems. We have also studied the effect that the different representations have on the index of the differential-algebraic system as {well as on} the associated Hamilton function. In general it can be seen that certain algebraic constraints do not contribute to the Hamiltonian and therefore can be separated from the system in an appropriate coordinate system. We have also characterized when  a general differential-algebraic system can be transformed to the different representations. { Several important tasks remain open. These include extensions to the case of non-regular   systems. These can be based on the results and methods in  Appendices A and B that are already proved for the non-square case. For systems with inputs and outputs, linear time-varying and nonlinear systems the extensions are currently  under consideration.}

\section*{Appendix}

\subsection*{Appendix A: Proof of Lemma~\ref{lem:PScanform}}
\proof
We present a proof in form of an algorithmic procedure that can be implemented as a numerical algorithm.

Step 1. Let $V_1$ and $W_1$ be real orthogonal matrices such that
\[ V_1^\top  P W_1=\mat{cc} \tilde P_{11} & 0 \\ 0 & 0 \rix, \ V_{1}^\top  S W= \mat{cc} \tilde S_{11} & \tilde S_{12} \\ \tilde S_{21} & \tilde S_{22} \rix,
\]
where $\tilde P_{11}\in \mathbb R^{\tilde m_1,\tilde m_1}$ is diagonal with positive diagonal elements.  This transformation can be constructed via a singular value decomposition of $P$ and a numerical rank decision.

Then, using the structure \eqref{lagrange}, it follows that $\tilde P_{11}^\top  \tilde S_{11}=\tilde S_{11}^\top  \tilde P_{11}$ and $\tilde S_{12}=0$.

Step 2. Let $\tilde V_2$, $\tilde W_2$ be real orthogonal matrices such that
\[
\tilde V_2^\top  \tilde S_{22} W_2= \mat{cc} \hat S_{22} & 0 \\  0 & 0\rix
\]
where $\hat S_{22}\in \mathbb R^{m_4,m_4}$ is diagonal with positive diagonal elements.  This transformation can be constructed again via a singular value decomposition  and a numerical rank decision. Set
\[
V_2=\mat{cc} I_{\tilde m_1} & 0 \\ 0 & \tilde V_2^\top  \rix, \ W_2 =\mat{cc} I_{\tilde m_1} & 0 \\ 0 & \tilde W_2\rix
\]
and form
\[
V_2^\top  V_1^\top  P W_1 W_2= \mat{ccc} \hat P_{11} & 0 & 0\\ 0 & 0 & 0 \\ 0 & 0 & 0 \rix,\
V_2^\top  V_1^\top  S W_1 W_2= \mat{ccc} \hat S_{11} & 0 & 0\\ \hat S_{21} & \hat S_{22} & 0 \\ \hat S_{31} & 0 & 0 \rix,
\]
with $\hat P_{11}^\top =\tilde P_{11}$.


Step 3. Let
\[
 V_3= \mat{ccc} \hat P_{11} & 0 & 0 \\ 0 & \hat S_{22}^{-\top} & 0 \\ 0 & 0 & I_{n-\tilde m_1-m_4}\rix,\
  W_3=\mat{ccc} I_{\tilde m_1} & 0 & 0 \\ -\hat S_{22}^{-1} \hat S_{21} & I_{m_4} & 0 \\ 0 & 0 & I_{m-\tilde m_1-m_4}\rix,
\]
and form
\[
V_3^{-1}V_2^\top  V_1^\top  P W_1 W_2 W_3= \mat{ccc} I_{\tilde m_1} & 0 & 0\\ 0 & 0 & 0 \\ 0 & 0 & 0 \rix,\
V_{3}^\top  V_2^\top  V_1^\top  S W_1 W_2 W_3= \mat{ccc} \check S_{11} & 0 & 0\\ 0 & I_{m_4} & 0 \\ \check S_{31} & 0 & 0 \rix,
\]
where by the structure \eqref{lagrange}  now $\check S_{11}$ is symmetric. Note that although we are working with nonorthogonal transformation matrices in this step the numerical errors can be controlled, since we are inverting diagonal matrices.

Step 4. Let
\[
\check V_{11}^\top \check S_{11} \check V_{11}=  \mat{ccc} I_{m_1} & 0 & 0 \\ 0 & -I_{m_2} & 0 \\ 0 & 0 & 0 \rix
\]
be the  canonical form of the skew-symmetric matrix $\check S_{11}$ under congruence which can be obtained by first computing the spectral decomposition and then scaling the nonsingular diagonal parts by congruence to be $\pm I$, see e.g.  \cite{LieM15a}.

Furthermore let
\[
\check S_{31} = \check V_4  \mat{ccc} S_{51} & S_{5,2} & S_{5,3}\\ 0 & 0 & 0\rix
\]
be a full rank decomposition partitioned accordingly, with $\check V_3$ real orthogonal. Then set
\[
 V_4= \mat{ccc} \check V_{11}^{-\top} & 0 & 0 \\ 0 & I_{m_4} & 0 \\ 0 & 0 & \check V_4\rix ,\
 W_4=\mat{ccc} \check V_{11} & 0 & 0 \\ & I_{m_4} & 0 \\ 0 & 0 & I_{m-\tilde m_1-m_4}\rix,
\]
and form
\begin{eqnarray*}
V_4^{-1} V_3^{-1}V_2^{-1} V_1^{-1} P W_1 W_2 W_3 W_4& =& \mat{ccccc} I_{m_1} & 0 & 0 & 0 & 0 \\ 0& I_{m_2} & 0 & 0 & 0 \\
0& 0 & I_{m_3} & 0 & 0 \\ 0 & 0 & 0 & 0 & 0 \\ 0 & 0 & 0 & 0 & 0 \rix,\\
V_4^{-1} V_3^{-1}V_2^{-1} V_1^{-1}  S W_1 W_2 W_3&=& \mat{cccccc}I_{m_1} & 0 & 0 & 0 & 0\\
0 & -I_{m_2} & 0 & 0 & 0\\
0 & 0 & 0 & 0 &  0\\
0 & 0 & 0 & I_{m_4} & 0 \\
S_{51} & S_{5,2} & S_{5,3} & 0  & 0\\
0 & 0 & 0 & 0 & 0
\rix,
\end{eqnarray*}
which is as claimed.
\eproof
\subsection*{Appendix B: Proof of Lemma~\ref{lem:KLcanform}}
\proof
The proof is similar to that of Lemma~\ref{lem:PScanform} just adapting to the different symmetries and transformation structure.
The following algorithmic procedure can be directly implemented as a numerical algorithm.

Step 1. Let $U_1$ and $V_1$ be real orthogonal matrices such that
\[ U_1^\top  K V_1=\mat{cc} \tilde K_{11} & 0 \\ 0 & 0 \rix, \ U_{1}^\top  L V_1^{-\top}= \mat{cc} \tilde L_{11} & \tilde L_{12} \\ \tilde L_{21} & \tilde L_{22} \rix,
\]
where $\tilde K_{11}\in \mathbb R^{\tilde \ell_1,\tilde \ell_1}$ is diagonal with positive diagonal elements.  This transformation can be constructed via a singular value decomposition of $K$ and a numerical rank decision.

Then, using the structure \eqref{dirac}, it follows that $\tilde K_{11}^\top  \tilde L_{11}=\tilde L_{11}^\top  \tilde K_{11}$ and $\tilde L_{21}=0$.

Step 2. Let $\tilde U_2$, $\tilde V_2$ be real orthogonal matrices such that
\[
\tilde U_2^\top  \tilde L_{22} V_2^{-\top}= \mat{cc} \hat L_{22} & 0 \\ 0 & 0 \rix,
\]
where $\hat L_{22}\in \mathbb R^{\ell_4,\ell_4}$ is diagonal with positive diagonal elements.  This transformation can be constructed again via a singular value decomposition  and a numerical rank decision. Set
\[
U_2=\mat{cc} I_{\tilde \ell_1} & 0 \\ 0 & \tilde U_2 \rix, \ V_2 =\mat{cc} I_{\tilde \ell_1} & 0 \\ 0 & \tilde V_2\rix
\]
and form
\[
U_2^\top  U_1^\top  K V_1 V_2= \mat{ccc} \hat K_{11} & 0 & 0\\ 0 & 0 & 0 \\ 0 & 0 & 0 \rix,\
U_2^\top  U_1^\top  L V_1^{-\top} V_2^{-\top}= \mat{ccc} \hat L_{11} & \hat L_{12} & \hat L_{13}\\ 0 & \hat L_{22} & 0\\ 0 & 0 & 0 \rix,
\]
with $\hat K_{11}=\tilde K_{11}$.

Step 3. Let
\[
 U_3= \mat{ccc} I_{\tilde \ell_1} & -\hat L_{22}^{-\top} \hat L_{12}^\top  & 0 \\ 0 & \hat L_{22}^{-\top} & 0 \\ 0 & 0 & I_{\ell-\tilde \ell_1-\ell_4}\rix,\
  V_3=\mat{ccc} \hat K_{11}^{-1} & 0 & 0 \\ 0 & I_{\ell_4} & 0 \\ 0 & 0 & I_{n-\tilde \ell_1-\ell_4}\rix,
\]
and form
\[
U_3^\top U_2^\top  U_1^\top  K V_1 V_2 V_3= \mat{ccc} I_{\tilde \ell_1} & 0 & 0\\ 0 & 0 & 0 \\ 0 & 0 & 0 \rix,\
U_{3}^\top  U_2^\top  U_1^\top  L V_1^{-\top} V_2^{-\top} V_3^{-\top}= \mat{ccc} \check L_{11} & 0 & \check L_{13}\\ 0 & I_{\ell_4} & 0 \\ \check 0 & 0 & 0 \rix,
\]
where by the structure \eqref{dirac}  now $\check L_{11}$ is skew-symmetric. Note that although we are working with nonorthogonal transformation matrices in this step, the numerical errors can be controlled since we are inverting diagonal matrices.

Step 4. Let

\[
\check U_{11}^\top \check L_{11} \check U_{11}=  \mat{ccc} 0 & I_{\ell_1} & 0 \\  -I_{\ell_1} & 0 &0\\ 0 & 0 & 0 \rix
\]
be the  canonical form of the skew-symmetric matrix $\check L_{11}$ under congruence which can be obtained by first computing the spectral decomposition and then scaling the nonsingular diagonal parts by congruence to be $\pm I$. This procedure is implemented in a numerically robust way in \cite{BenBFMW00}.
Furthermore let
\[
\check L_{13} =   \mat{cc} L_{15} & 0 \\ L_{25} & 0\\ L_{35} & 0 \rix \check V_3
\]
be a full rank decomposition partitioned accordingly, with $\check V_3$ real orthogonal. Then set
\[
 U_4= \mat{ccc} \check U_{11} & 0 & 0 \\ 0 & I_{\ell_4} & 0 \\ 0 & 0 & I_{\ell-\tilde \ell_1 -\ell_4}\rix,\
 V_4=\mat{ccc} \check U_{11}^{-\top} & 0 & 0 \\ & I_{\ell_4} & 0 \\ 0 & 0 & \check V_3\rix,
\]
and form
\begin{eqnarray*}
U_4^{T} U_3^{T} U_2^\top  U_1^\top  K V_1 V_2 V_3 V_4&=& \mat{ccccc} I_{\ell_1} & 0 & 0 & 0 & 0 \\ 0& I_{\ell_1} & 0 & 0 & 0 \\
0& 0 & I_{\ell_3} & 0 & 0 \\ 0 & 0 & 0 & 0 & 0 \\ 0 & 0 & 0 & 0 & 0 \rix,\\
U_4^{T} U_3^{T} U_2^\top  U_1^\top  L V_1^{-\top} V_2^{-\top} V_3^{-\top} V_4^{-\top}&=& \mat{cccccc}0 & I_{\ell_1} & 0 & 0 & L_{15} & 0\\
-I_{\ell_1} & 0 & 0& 0 & L_{25} & 0\\
0 & 0 & 0 & 0 &  L_{35} & 0\\
0 & 0 & 0 & I_{\ell_4} & 0 & 0\\
0 & 0 & 0 & 0 & 0 & 0
\rix,
\end{eqnarray*}
which is as claimed. %
\eproof

\subsection*{Appendix C: Proof of Theorem~\ref{thm:sinind2}}
\proof
 We again present a proof that can be implemented as a numerical algorithm.

Step 1. Let $U_1$ and $V_1$ be real orthogonal matrices such that
\[ U_1^\top  E V_1=\mat{cc} \tilde E_{11} & 0 \\ 0 & 0 \rix, \ U_{1}^\top  A V_1= \mat{cc} \tilde A_{11} & \tilde A_{12} \\ \tilde A_{21} & \tilde A_{22} \rix,
\]
where $\tilde E_{11}\in \mathbb R^{\tilde n_1,\tilde n_1}$ is diagonal with positive diagonal elements.  This transformation can be constructed via a singular value decomposition of $E$ and a numerical rank decision.

Step 2. Let $\tilde U_2$, $\tilde V_2$ be real orthogonal matrices such that
\[
\tilde U_2^\top  \tilde A_{22} \tilde V_2= \mat{cc} \hat A_{22} & 0 \\ 0 & 0 \rix,
\]
where $\hat A_{22}\in \mathbb R^{n_3,n_3}$ is diagonal with positive diagonal elements.  This transformation can be constructed again via a singular value decomposition  and a numerical rank decision. Set
\[
U_2=\mat{cc} I_{\tilde n_1} & 0 \\ 0 & \tilde U_2 \rix, \ V_2 =\mat{cc} I_{\tilde n_1} & 0 \\ 0 & \tilde V_2\rix
\]
and form
\[
U_2^\top  U_1^\top  E V_1 V_2= \mat{ccc} \hat E_{11} & 0 & 0\\ 0 & 0 & 0 \\ 0 & 0 & 0 \rix,\
U_2^\top  U_1^\top  A V_1^{-\top} V_2= \mat{ccc} \hat A_{11} & \hat A_{12} & \hat A_{13}\\ A_{21} & \hat A_{22} & 0\\ A_{31}& 0 & 0 \rix,
\]
with $\hat E_{11}=\tilde E_{11}$.

The regularity of the pencil implies that $A_{13}$ has full column rank $n_1=\tilde n_1-n_3$ and that $A_{31}$ has full row rank $n_1=\tilde n_1-n_3$, because otherwise there would be common right or left nullspace, respectively.

Step 3.
Let $\tilde U_{31}$, $\tilde V_{31}$, $\tilde U_{13}$, $\tilde V_{13}$ be real orthogonal matrices of appropriate dimensions such that
\[
\tilde U_{13}^\top  \tilde A_{13} \tilde V_{13}= \mat{c} \hat A_{14}  \\ 0  \rix,\
\tilde U_{31}^\top  \tilde A_{31} \tilde V_{31}= \mat{cc} \hat A_{41} & 0  \rix
\]
where $\hat A_{41}\in \mathbb R^{n_1,n_1}$  and $\hat A_{14}\in \mathbb R^{n_1,n_1}$ are diagonal with positive diagonal elements.  This transformation can be constructed again via a singular value decomposition  and  numerical rank decisions. Set
\[
U_3=\mat{ccc} \tilde U_{13} & 0 &0 \\ 0 & I_{\ell_3} & 0  \\ 0 & 0 & \tilde U_{31}\rix, \
V_3=\mat{ccc} \tilde V_{31} & 0 &0 \\ 0 & I_{\ell_3} & 0  \\ 0 & 0 & \tilde V_{13}\rix,
\]
then $U_3^\top U_2^\top  U_1^\top  E V_1 V_2 V_3$ and $U_3^\top U_2^\top  U_1^\top  A V_1 V_2 V_3$ are as claimed in \eqref{cformind2}. The invertibility of $E_{22}$ then follows from the assumption that the pencil has index at most two, see \cite{KunM06}.
\eproof

\end{document}